%% file: metasurfaces-ot.tex
\documentclass[11pt,a4]{amsart}
\usepackage{amsmath,amsfonts,latexsym,graphicx,amssymb,url}
\usepackage{hyperref,pdfsync}
\usepackage{amsmath}
\usepackage[active]{srcltx}
\usepackage{cases,soul}
\usepackage{wasysym, enumerate}
\usepackage{lscape,color,datetime,subfigure,graphicx,wrapfig}
\usepackage{cleveref,float}
\usepackage{cancel}

\newcommand{\dist}{\text{dist}} 

\newcommand{\diver}{\text{\rm div}}
\newcommand{\R}{{\mathbb R}} 

\newcommand{\N}{{\mathbb N}}

\renewcommand{\(}{\left(}
\renewcommand{\)}{\right)}
\newcommand{\e}{{\bf e}}
\newcommand{\m}{{\bf m}}
\newcommand{\x}{{\bf x}}
\renewcommand{\b}{{\bf b}}
\newcommand{\g}{{\bf g}}
\newcommand{\Lag}{\text{Lag}}
\newcommand{\Pow}{\text{Pow}}
\newcommand{\proj}{\text{P}}
\newcommand{\eps}{{\varepsilon}}
\renewcommand{\epsilon}{\varepsilon}
\newcommand{\nr}[1]{\vert#1\vert}
\newcommand{\X}{{\mathcal X}}
\newcommand{\Y}{{\mathcal Y}}

\newtheorem{theorem}{Theorem}[section]

\newtheorem{corollary}[theorem]{Corollary}
\newtheorem{lemma}[theorem]{Lemma}
\newtheorem{proposition}[theorem]{Proposition}
\newtheorem{definition}[theorem]{Definition}
\newtheorem{remark}[theorem]{Remark}

\numberwithin{equation}{section}

\title[Metasurfaces and Optimal transport, ~~~\today]{Metasurfaces and Optimal transport}
\author[Guti\'errez, Huang, M\'erigot, and Thibert]{Cristian E. Guti\'errez, Qingbo Huang, \\Quentin M\'erigot and Boris Thibert}
\thanks{\today}
\address{Department of Mathematics\\Temple University\\Philadelphia, PA 19122}
\email{gutierre@temple.edu}
\address{Department of Mathematics and Statistics\\Wright State University\\Dayton, OH 45435}
\email{qingbo.huang@wright.edu}
\address{Universit\'e Paris-Saclay, CNRS, Laboratoire de math\'ematiques d'Orsay, 91405, Orsay, France \and
Institut universitaire de France (IUF)}
\email{quentin.merigot@universite-paris-saclay.fr}
\address{Universit\'e  Grenoble Alpes, CNRS, Grenoble INP, LJK, 38000 Grenoble, France}
\email{Boris.Thibert@univ-grenoble-alpes.fr}

\begin{document}
\begin{abstract}
This paper provides a theoretical and numerical approach to show existence, uniqueness, and the numerical determination of metalenses refracting radiation with energy patterns. 
The theoretical part uses ideas from optimal transport and for the numerical solution we study and implement a damped Newton algorithm to solve the semi discrete problem. 
A detailed analysis is carried out to solve the near field one source refraction problem and extensions to the far field are also mentioned.
\end{abstract}
\maketitle

\input{section1-intro}
\input{section2-theory}

\input{section3-mapping}
\input{section4-numerics}
\input{section5-farfield}

\end{document}

%% file: section1-intro.tex
\section{Introduction}
Metalenses or metasurfaces are ultra thin surfaces with arrangements of  
nano scattering structures designed to focus light in imaging.
They introduce
abrupt phase changes over the scale of the wavelength along the
optical path to bend light in unusual ways.  
This is in
contrast with conventional lenses, where the question is to determine its faces so that a gradual change of phase accumulates as the wave
propagates inside the lens, reshaping the scattered wave at will.
These nano structures are engineered by adjusting their shape, size, position and orientation, and arranged on the surface (typically a plane) in the form of  
tiny pillars, rings, and others dispositions, working together to manipulate
light waves as they pass through. 
The subject of metalenses is an important area of current
research, one of the nine runners-up for Science's Breakthrough of the Year
2016 \cite{science-runner-ups-2016}, and is potentially useful in imaging applications.
Metalenses are thinner than a sheet of paper and far lighter than
glass, and they could revolutionize optical imaging devices from microscopes to virtual reality
displays and cameras, including the ones in smartphones; see  
\cite{science-runner-ups-2016}, \cite{metalenses-cameras}, and \cite{chen-capasso-21}.

Mathematically, a metalens can be described as pair $(\Gamma,\phi)$, where $\Gamma$ is a surface in 3-d space given by the graph of a $C^2$ function $u$, and $\phi$ is a function, called the {\it phase discontinuity}, defined in a small neighborhood of $\Gamma$. 
The knowledge of $\phi$ yields the kind of
arrangements of the nano structures on the surface that are needed for a specific refraction job. 
Refraction here acts following  the generalized Snell law \eqref{eq:generalized snell law}. 
For example, if $\Gamma=\Pi$ is a plane and the phase $\phi$ has the form  \eqref{eq:definition of the cost cXY}, then the metasurface $(\Gamma,\phi)$ refracts all rays from the origin $O$ into the point $Y$.

The question considered in this paper concerns existence, uniqueness and the numerical determination of metalenses refracting radiation with energy patterns.
We state the problem in the near field case; the far field is explained in Section \ref{sec:refraction into the far field}. 
Precisely, suppose radiation is emanating from a point source $O$, below a given surface $\Gamma$, with intensity $f(x)$ for each $x\in \Omega$ a domain of the unit sphere $S^2$. Furthermore, $T$ is a compact set, above the surface $\Gamma$, and a distribution of energy on $T$ is given by a Radon measure $\mu$ so that $\int_\Omega f(x)\,dx=\mu(T)$. We denote $\mathcal T(E)$ (see Definition \ref{def:refractor mapping})  the collection of points in $\Omega$ that are refracted into $E$ in accordance with the generalized Snell law \eqref{eq:generalized snell law}. Then, under what circumstances is there a phase discontinuity function $\phi$ defined in a neighborhood of $\Gamma$ so that the metalens $(\Gamma,\phi)$ refracts $\Omega$ into $T$ and satisfies the energy conservation 
\[
\int_{\mathcal T (E)}f(x)\,dx=\mu(E)
\]
for each $E\subset T$? 
We will solve this problem using ideas from
optimal transport.  Our first result is existence and uniqueness of
solutions in the semi-discrete case, that is, when $\mu$ is a finite
combination of delta functions and $\Gamma$ is a plane, Theorem
\ref{thm:existence of metasurfaces near field}.  A relative visibility
condition between $\Omega$ and $T$ is needed, condition
\eqref{eq:condition for single valued}, and to obtain our results we
use \cite{gutierrez-huang:nearfieldrefractor}.  We also provide a
numerical solution of the semi-discrete problem using a damped Newton
algorithm, introduced in Section \ref{sec:algorithm}.  This requires a
careful analysis of the refractor mapping, Definition
\ref{def:refractor mapping}, and the Laguerre cells in
\eqref{eq:definition of Laguerre cell}.

We mention that the phase discontinuity functions needed to design metalenses for various refraction and reflection problems with prescribed distributions of energy satisfy partial differential equations of Monge-Amp\`ere type which are derived and 
studied in \cite{gutierrez-pallucchini:metasurfacesandMAequations}. 
Equations of these type also appear naturally in solving problems involving aspherical lenses, see \cite{gutierrez-huang:reshaping light beam}, \cite{abedin-gutierrez-tralli:regularity-refractors},  \cite{gutierrez-sabra:asphericallensdesignandimaging}, \cite{gutierrez-sabra:freeformgeneralfields}, and references therein.

The paper is organized as follows.
 Section \ref{sec:existence and uniqueness} contains a precise description of the problem and the existence and uniqueness results. The Laguerre cells for our problem and the analysis of the refractor mapping is the contents of Section \ref{sec:kantorovitch formulation}. To handle the singular set  \eqref{eq:singular set} where the Laguerre cells intersect the boundary of $\Omega$, we assume that the target $T$ is contained in a plane parallel to the metasurface and the boundary $\partial \Omega$ is not a conic section, see Remark \ref{rmk:assumptions on Omega}.
In Section \ref{sec:computation of
  the Laguerre cells} we show that 2D Laguerre cells, which are complicated objects,
can be computed in a simpler way from a 3D power diagram, a tessellation of the 3D space
into convex polyhedra. This leads to an effective method to solve the
near-field refractor problem, which is tested on a few cases.
Finally, Section \ref{sec:refraction into the far field} contains the
far field case for both collimated and point sources.

\subsection*{Acknowledgements}
A large part of this work was carried out while the first named author was visiting the Universities of Paris-Sud and Grenoble Alpes during Fall 2019 on sabbatical from Temple University under NSF grant DMS-1600578.
He would like to warmly thank his co-authors Quentin M\'erigot and Boris Thibert and their institutions for the hospitality and support during his visit. 
The last two named authors acknowledge support from the French Agence National de la Recherche through the project MAGA (ANR-16-CE40-0014).

%% file: section2-theory.tex
\section{Refraction from one point into a near field target}
\label{sec:existence and uniqueness}

\subsection{Generalized Snell's law}
\label{sec:preliminaries}

Let $S\subseteq \R^3$ be a smooth surface defined implicitely by the equation
$\psi(x)=0$, where $x=(x_1,x_2,x_3)\in \R^3$, and let $\phi$ be a
function defined in a small neighborhood of $S$.  The region below $S$
is made up of an homogeneous material \textsc{I} with refractive index
$n_1$, and the region above $S$ is made up of an homogeneous
material \textsc{II} having refractive index $n_2$. From a point $A$
in \textsc{I}, a wave is emitted, it strikes $S$ at some point $X$,
and is then transmitted to a point $B$ in medium \textsc{II}.  Let us
fix $A$ and $B$ and we want to minimize
\[
n_1\,|x-A|+n_2\,|x-B|-\phi(x)
\]
over all $x\in S$, i.e., with $\psi(x)=0$.  From the existence of
Lagrange multipliers, the gradient vector
$\nabla_x\(n_1\,|x-A|+n_2\,|x-B|-\phi(x)\)$ must be parallel to the
normal vector $\nabla \psi(x)$ at all critical points $x$.  That is,
the vector product
$\nabla_x\(n_1\,|x-A|+n_2\,|x-B|-\phi(x)\)\times \nu(x)=0$, with
$\nu(x)$ being the normal to $S$. If we set $\x=\dfrac{x-A}{|x-A|}$,
and $\m=\dfrac{B-x}{|x-B|}$, then $\(n_1\,\x- n_2\,\m
-\nabla \phi(x)\)\times \nu(x)=0$.  This means the vectors are
multiple one from the other, that is, we obtain the {\it generalized
Snell law}
\begin{equation}\label{eq:generalized snell law}
n_1\,\x- n_2\,\m =\lambda\,\nu(x)+\nabla \phi(x),
\end{equation}
for some $\lambda\in \R$; the function $\phi$ is called the phase discontinuity; for a derivation of this law using wave fronts see \cite{gps}.

\subsection{Formulation of the problem}
Let $\Pi$ denote the plane $x_3=\alpha$ in $\R^3$, $\alpha>0$.
Rays emanate from the origin $O$ with directions $x\in \Omega_0\subset S^2$, a compact domain, and intensity $f(x)$.
Given $x\in \Omega_0$, let $X\in \Pi$ so that $x=X/|X|$, and  
set 
\[
\Omega=\{X\in \Pi: X/|X|\in \Omega_0\},
\] 
establishing a one to one correspondence between $\Omega_0\subset S^2$ and the compact domain $\Omega\subset \Pi$.
From \cite[Section 7.A]{gps}, given a point $Y$ above the plane $\Pi$, the phase discontinuity $\phi$ so that the metasurface $(\Pi,\phi)$ refracts rays from $O$ into $Y$, with $\phi$ tangential to $\Pi$, i.e., $\nabla \phi(X)\cdot (0,0,1)=0$ \footnote{Each $\phi(X)=|X|+|X-Y|+h(x_3)$ satisfies \eqref{eq:generalized snell law} with $n_1=n_2=1$ but is not tangential to $\Pi$ unless $h$ is constant.}, is given by 
\begin{equation}\label{eq:definition of the cost cXY}
\phi(X)=|X|+|X-Y|:=c(X,Y),\quad X\in \Pi.
\end{equation}
Let $T$ be a compact domain in $\R^3$ above the plane $\Pi$, with $\dist(T,\Pi)>0$, $T$ is referred as the target or receiver.  
\begin{definition}[Admissible phase near field]\label{def:admissible phase near field}
The function $\phi:\Omega\to \R$ is an admissible phase refracting $\Omega$ into $T$ if for each $X_0\in \Omega$ there exists $b\in \R$ and $Y\in T$ such that 
\[
\phi(X)\leq c(X,Y)+b \quad \forall X\in \Omega, \quad \phi(X_0)= c(X_0,Y)+b.
\]
In this case, we say that $c(X,Y)+b$ supports $\phi$ at $X_0$.
\end{definition}
Is easy to see that admissible phases are Lipschitz continuous, $|\phi(X)-\phi(Y)|\leq 2\,|X-Y|$ for all $X,Y\in \Omega$.

\begin{definition}[Refractor mapping]\label{def:refractor mapping}
If $\phi:\Omega\to \R$ is an admissible phase, then for each $X_0\in \Omega$ we define the set valued mapping 
\[
\mathcal N_\phi(X_0)=\{Y\in T: \text{there exists $b\in \R$ such that $c(X,Y)+b$ supports $\phi$ at $X_0$}\},
\]
and for each $Y\in T$ we set
\[
\mathcal T_\phi(Y)=\mathcal N_\phi^{-1}(Y)=\{X\in \Omega: Y\in \mathcal N_\phi(X)\}.
\]
$\mathcal N_\phi$ is the refractor mapping and $\mathcal T_\phi$ the tracing mapping.
If $E\subset \Omega$, then 
\[
\mathcal T_\phi(E)=\cup_{Y\in E}\mathcal T_\phi(Y)
=\{X\in \Omega:\mathcal N_\phi(x)\cap E\neq \emptyset\}.
\]
\end{definition}

Now, let $f\in L^1(\Omega_0)$ with $f>0$ a.e., and set $\rho(X)$ to be the density induced by $f$ for $X\in \Omega$; $\rho(X)=|X|^{-3}\,f(X/|X|)\,X\cdot e$ with $e=(0,0,1)$.
In addition, let $\mu$ be a Radon measure in $T$ satisfying the energy conservation condition 
\begin{equation}\label{eq:energy conservation condition}
\int_\Omega \rho(X)\,dX=\mu(T).
\end{equation}
{\it The problem we consider here is that of finding an admissible phase function $\phi:\Omega\to \R$ so that 
\begin{equation}\label{eq:energy conservation condition local}
\int_{\mathcal T_\phi(E)} \rho(X)\,dX=\mu(E)
\end{equation}
for each Borel set $E\subset T$. This means the metasurface $(\Omega,\phi)$ refracts $\Omega$ into $T$ satisfying the energy conservation condition \eqref{eq:energy conservation condition local}.
In the following section we prove existence of solutions to this problem.}

In order to get existence and uniqueness of solutions, we assume  that $|\partial \Omega|=0$\footnote{Recall that $\Omega$ is on the plane $\Pi$, $\partial \Omega$ denotes its  boundary in $\Pi$ and $|\partial \Omega|$ denotes its two dimensional Lebesgue measure.} 
and that the points $Y_1,\hdots,Y_N$ satisfy the following condition, which holds for instance if the target $T$ is contained in the plane $x_3=\beta$ with $\beta>\alpha$:
\begin{equation}\label{eq:condition for single valued}
\forall (X,Y_1,Y_2) \in \Omega\times T\times T,~~ Y_1 \neq Y_2 \Longrightarrow  X,Y_1,Y_2 \hbox{ are not aligned.}
\end{equation}
It will be proved in Section \ref{sec:existence of solutions} below that \eqref{eq:condition for single valued} implies that $\mathcal N_\phi$ is single valued for almost every $X$. 
  


We then have the following existence and uniqueness theorem.

%
 
\begin{theorem}\label{thm:existence of metasurfaces near field}
Let $\Omega$ be a compact connected domain on the plane $\Pi=\{x_3=\alpha\}$, $\alpha>0$, with $|\partial \Omega|=0$, and let $Y_1,\cdots , Y_N$ be distinct points in the target $T$ laying above $\Pi$ and satisfying \eqref{eq:condition for single valued} with $\dist(\Omega,T)>0$.
Further, let $g_{1},\cdots , g_{N}$ be positive numbers, and $\rho \in L^{1}(\Omega)$ satisfying \eqref{eq:energy conservation condition} with the measure $\mu=\sum_{i=1}^N g_i\,\delta_{Y_i}$.

Then 
given any $b_1\in \R$, there exist unique numbers $b_{2},\cdots , b_{N}$ such that the function
\begin{equation}\label{eq:solutionsecondbdryMA near field}
\phi(X)=\min_{1\leq i \leq N}\{ |X|+|X-Y_i|+b_i\}
\end{equation}
solves \eqref{eq:energy conservation condition local}, that is, the metasurface $(\Omega,\phi)$ refracts $\Omega$ into $T$. 
\end{theorem}

\begin{remark}[Relation to optimal transport] One can check that the
the couple $(\phi,\b)$ constructed in Theorem~\ref{thm:existence of metasurfaces near field} is the solution to the following maximization problem:
$$ \max \left\{ \int \phi d\rho + \sum_{1\leq i\leq N} b_i g_i : \phi\in C^0(\Omega), \b\in\R^N, \forall X\in T, \forall i\in \{1,\hdots,N\}, \phi(X) + b_i \leq c(X,Y_i) \right\},$$
where $c(X,Y_i) = |X|+|X-Y_i|$.
Thus, solving \eqref{eq:energy
conservation condition local} amounts to solving the Kantorovich dual
of the optimal transport problem
$$ \min  \left\{ \int c(X,Y) d \gamma(X,Y) : \gamma \in \Pi(\rho,\mu) \right\}. $$
where $\Pi(\rho,\mu)$ denotes the set of \emph{transport plans}
between $\rho$ and $\mu$, i.e. probability measures on $\Omega\times
T$ with respective marginals $\rho$ and $\mu$.
\end{remark}

\subsection{Existence of solutions}\label{sec:existence of solutions}
We will use the results from \cite[Section 2]{gutierrez-huang:nearfieldrefractor}, and we 
first recall some notation from there. Suppose $\X, \Y$ are compact metric spaces and $\omega$ is a Radon measure in $\X$. $\Y^\X$ denotes the class of all set valued mappings $\Phi:\X\to \mathcal P(\Y)$ that are single valued for almost all points in $\X$ with respect to the measure $\omega$.
We say that $\Phi\in \Y^\X$ is continuous at the point $x_0\in \X$ if given $x_k\to x_0$ and $y_k\in \Phi(x_k)$ there exists a subsequence $y_{k_j}$ and $y_0\in \Phi(x_0)$ such that $y_{k_j}\to y_0$ as $j\to \infty$. We also denote
\begin{align*}
C(\X,\Y)=\{\Phi\in \Y^\X: \text{$\Phi$ is continuous in $\X$}\},
\text{ and }
C_s(\X,\Y)=\{\Phi\in C(\X,\Y):\Phi(\X)=\Y\}.
\end{align*}

With this set up, we let $\X=\Omega$, $\Y=T$, $\omega=\rho\,dx$, and 
to solve our problem we introduce the class
\[
\mathcal F=\{\phi:\Omega\to \R:\text{$\phi$ is an admissible phase refracting $\Omega$ into $T$}\}.
\]
From the definitions introduced above, it is easy to see that $\mathcal F$ satisfies the following properties like the ones introduced in \cite[Sections 2.1 and 2.3]{gutierrez-huang:nearfieldrefractor} (with the same labeling so the reader can compare):
\begin{enumerate}
\item[(A1)] if $\phi_1,\phi_2\in \mathcal F$, then $\phi_1\wedge\phi_2=\min\{\phi_1,\phi_2\}\in \mathcal F$,
\item[(A2)] if $\phi_1(x_0)\leq \phi_2(x_0)$, then $\mathcal N_{\phi_1}(x_0)\subset \mathcal N_{\phi_1\wedge \phi_2}(x_0)$,
\item[(A3'')] for each $Y\in T$ and each $b\in (-\infty,\infty)$, the functions $c(\cdot,Y)+b \in \mathcal F$ satisfy the following
\begin{enumerate}
\item[(a)] $Y\in \mathcal N_{c(\cdot,Y)+b}(X)$ for all $X\in \Omega$,
\item[(b)] $c(\cdot,Y)+b\leq c(\cdot,Y)+b'$ for all $b\leq b'$,
\item[(c)] for each $Y\in T$, $c(X,Y)+b\to +\infty$ uniformly for $X\in \Omega$ as $b\to +\infty$,
\item[(d)] for each $Y\in T$, $\max_{X\in \Omega} |c(X,Y)+b-\(c(X,Y)+b'\)|\to 0$ as $b'\to b$.
\end{enumerate}
\end{enumerate}

In addition, we need to verify that 
\begin{equation}\label{eq:N phi is in Cs}
\mathcal N_\phi\in C_s(\Omega,T)
\end{equation}
for each $\phi\in \mathcal F$. We first verify that $\mathcal N_\phi\in T^\Omega$, that is, $\mathcal N_\phi$ is single valued for almost all $X\in \Omega$ with respect to Lebesgue measure. Indeed, suppose that $\mathcal N_\phi(X_0)$ contains more than one point, say $Y_1,Y_2\in \mathcal N_\phi(X_0)$ with $Y_1\neq Y_2$ and $X_0\in \mathring \Omega$. Will show that $X_0$ is a singular point of $\phi$.
We have that for some $b_1,b_2\in \R$, $\phi(X)\leq |X|+|X-Y_i|+b_i$ for all $X\in \Omega$ with equality at $X=X_0$, $i=1,2$.
Since $\dist(\Omega,O)>0$ and $\dist(\Omega,T)>0$, the support functions $|X|+|X-Y_i|$ are both smooth in the variable $X=(x_1,x_2,\alpha)$. If $X_0$ were not a singular point of $\phi$, then $\phi$ has a tangent plane at $X_0$ which must coincide with both tangent planes to $|X|+|X-Y_i|+b_i$, $i=1,2$, at $X_0$.
So $\left. \nabla_{x_1,x_2}\(|X|+|X-Y_1|\)\right|_{X=X_0}=\left. \nabla_{x_1,x_2}\(|X|+|X-Y_2|\)\right|_{X=X_0}$ obtaining $\dfrac{X_0-Y_1}{|X_0-Y_1|}=\dfrac{X_0-Y_2}{|X_0-Y_2|}$ which implies that the vectors $X_0-Y_1$ and $X_0-Y_2$ are multiples one of the other. This contradicts \eqref{eq:condition for single valued} showing that $X_0$ is a singular point to $\phi$. Since $\phi$ is Lipschitz, the set of points where $\phi$ is not differentiable has measure zero and therefore \eqref{eq:condition for single valued} implies that $\mathcal N_\phi$ is single valued for almost all $X\in \Omega$. 

Second, we verify that $\mathcal N_\phi$ is continuous in the sense given at the beginning of this section.
Indeed, let $X_k\in \Omega$ with $X_k\to X_0$, and let $Y_k\in \mathcal N_\phi(X_k)$. We then have $\phi(X)\leq c(X,Y_k)+b_k$ for all $X\in \Omega$ for some $b_k$ with equality at $X=X_k$. 
Since $\Omega,T$ are compact, $\phi$ is continuous in $\Omega$, and $c$ is continuous on $\Omega\times T$, selecting subsequences is easy to obtain the desired continuity. 
We also have that $\mathcal N_\phi$ is continuous at each $\phi\in \mathcal F$, i.e., if $\phi_k\to \phi$ uniformly in $\Omega$, $\phi_k\in \mathcal F$, $X_0\in \Omega$ and $Y_k\in \mathcal N_{\phi_k}(X_0)$, then there exists a subsequence $Y_{k_j}\to Y_0\in T$ with $Y_0\in \mathcal N_\phi(X_0)$. In fact, there exists $b_k\in \R$ such that $\phi_k(X)\leq c(X,Y_k)+b_k$ for all $X\in \Omega$ with equality at $X=X_0$. That is, $b_k=\phi_k(X_0)-c(X_0,Y_k)=\phi_k(X_0)-\phi(X_0)+\phi(X_0)-c(X_0,Y_k)$, a quantity uniformly bounded in $k$ from the uniform convergence and since $\Omega,T$ are compact. Since $T$ is compact, there is a subsequence $Y_{k_j}\to Y_0$ for some $Y_0\in T$, and so there is a subsequence $b_{{k_j}_\ell}\to b_0$ as $\ell\to \infty$ obtaining that $c(X,Y_0)+b_0$ supports $\phi$ at $X_0$, that is, $Y_0\in    N_\phi(X_0)$.

To continue verifying \eqref{eq:N phi is in Cs}, we next show that $\mathcal N_\phi(\Omega)=T$.
Let $Y\in T$. Since $\phi$ is continuous in $\Omega$ we have 
\[
\phi(X)\leq \max_\Omega \phi\leq |X|+|X-Y|+b
\]
for all $X\in \Omega$ and for $b\to \infty$. 
Choose 
\[
b_0=\inf \{b:\phi(X) \leq |X|+|X-Y|+b,\forall X\in \Omega\}.
\]
Then $\phi(X_0)=|X_0|+|X_0-Y|+b_0$ for some $X_0\in \Omega$, that is, $\mathcal N_\phi(X_0)=Y$.

To prove existence of solutions to the problem \eqref{eq:energy conservation condition local} 
we proceed as follows.
We recall that in \cite[Sect. 2.3]{gutierrez-huang:nearfieldrefractor} we have considered the case convex infinity. However, to show existence in the current case, we need to consider the case concave infinity which is not explicitly written in \cite{gutierrez-huang:nearfieldrefractor} but it follows along similar lines as the convex infinity case. Indeed, to obtain existence of solutions we argue as in the proof of \cite[Theorem 2.12]{gutierrez-huang:nearfieldrefractor} using the family $\mathcal F^*=\{e^{-\phi},\phi\in \mathcal F\}$ together with the mapping $\mathcal N^*_{e^{-\phi}}=\mathcal N_\phi$
which now converts the case concave-infinity into the convex case from \cite[Sect. 2.2]{gutierrez-huang:nearfieldrefractor} with the family of supporting functions $e^{-c(X,Y)-b}$ for $b\in (-\infty,\infty)$. Hence existence in our case will follow applying \cite[Theorem 2.9]{gutierrez-huang:nearfieldrefractor} with the family $\mathcal F^*$ if we are under its hypotheses.
That is, we need to show that there exists $(b_1^0,b_2^0,\cdots ,b_N^0)$ such that
\[
e^{-c(X,Y_1)-b_1^0}\leq \min_{2\leq i\leq N} e^{-c(X,Y_i)-b_i^0}
\]
which is equivalent to
\[
\max_{2\leq i\leq N} c(X,Y_i)+b_i^0\leq c(X,Y_1)+b_1^0.
\]
Indeed, if $b_1^0$ is arbitrarily chosen then by continuity we can pick $b_2^0,\cdots ,b_N^0\to -\infty$ such that the last inequality holds.
 Therefore,  \cite[Theorem 2.9]{gutierrez-huang:nearfieldrefractor} implies the existence part of Theorem~\ref{thm:existence of metasurfaces near field}.

\subsection{Uniqueness of solutions to Theorem \ref{thm:existence of metasurfaces near field}}\label{sec:uniqueness of solutions}
We begin with a lemma.
\begin{lemma}\label{lm:support holds on non singular points}
Suppose \eqref{eq:condition for single valued} holds and $|\partial \Omega|=0$.
Let $\phi(X)=\min_{1\leq i\leq N}c(X,Y_i)+b_i$ with $Y_i\in T$ distinct points. 
If $X\in \mathcal T_{\phi}(Y_j)$, then $c(\cdot,Y_j)+b_j$ supports $\phi$ at $X$ or $X$ belongs to the set where $\mathcal N_{\phi}$ is not single valued (a set of measure zero). 
\end{lemma}
\begin{proof}
Let $X\in \mathcal T_{\phi}(Y_j)$, then there exists $\bar b$ such that $c(\cdot,Y_j)+\bar b$ supports $\phi$ at $X$. 
Then $c(X,Y_j)+\bar b=\phi(X)\leq c(X,Y_k)+b_k$ for all $k$, and in particular for $k=j$, so we get $\bar b\leq b_j$.
If $\bar b=b_j$, then $c(\cdot,Y_j)+b_j$ supports $\phi$ at $X$. 
If $\bar b<b_j$, then $\phi(Z)\leq c(Z,Y_j)+\bar b<c(Z,Y_j)+b_j$ for all $Z\in \Omega$ and hence $\phi(Z)=\min_{k\neq j}c(Z,Y_k)+b_k$ for all $Z\in \Omega$. 
In particular, when $Z=X$, there is $k\neq j$ such that $\phi(X)=c(X,Y_k)+b_k$. 
This means $c(\cdot,Y_k)+b_k$ supports $\phi$ at $X$, so $X\in \mathcal T_{\phi}(Y_k)$.
Hence $Y_j,Y_k\in \mathcal N_{\phi}(X)$ and so $X$ belongs to the set where $\mathcal N_{\phi}$ is not single valued which is a set of measure zero.
\end{proof}

For $\phi\in \mathcal F$, let $\mathcal Q=\{x\in \Omega:\mathcal N_{\phi}(x)\text{ is not singleton}\}$ which has measure zero because \eqref{eq:condition for single valued} holds and $|\partial \Omega|=0$. 
We also have that $\mathcal T_{\phi}(K)$ is compact for each $K$ compact subset of $\Omega$. 

\begin{lemma}\label{lm:boundary of tracing contained in singular set}
Let $\phi(X)=\min_{1\leq i\leq N}c(X,Y_i)+b_i$ with $Y_i\in T$ distinct points.
Then for each $F\subsetneqq \{Y_1,\cdots ,Y_N\}$, $F\neq \emptyset$ we have
\[
\partial \(\mathcal T_{\phi}(F)\)\subset \mathcal Q,
\]
where $\partial \(\cdot \)$ denotes the boundary.
\end{lemma}
\begin{proof}
Let $Z\in \partial \(\mathcal T_{\phi}(F)\)$. Then for each $m\geq 1$ the open ball $B_{1/m}(Z)$ satisfies 
$B_{1/m}(Z)\cap \(\mathcal T_{\phi}(F)\)^c\neq \emptyset$. 
Since $\mathcal T_{\phi}(F)$ is closed, $ \(\mathcal T_{\phi}(F)\)^c$ is open and so $B_{1/m}(Z)\cap \(\mathcal T_{\phi}(F)\)^c$ is a non empty relatively open set and so it has positive measure for all $m$.
Since $|\mathcal Q|=0$, it follows that the set $B_{1/m}(Z)\cap \(\mathcal T_{\phi}(F)\)^c\cap \mathcal Q^c$
has positive measure and we pick $Z_m$ in that set.
Hence $\mathcal N_\phi$ is single valued at $Z_m$, and $\mathcal N_\phi(Z_m)\cap F=\emptyset$.
Then it follows that $\mathcal N_\phi(Z_m)$ can take only values in $F^c$, and since $\{Y_1,\cdots ,Y_N\}$ is a finite set
there is subsequence such that $\mathcal N_\phi\(Z_{m_j}\)=Y_\ell$ for some $Y_\ell\in \{Y_1,\cdots ,Y_N\}\setminus F$. 
From Lemma \ref{lm:support holds on non singular points}, this means that $\phi(X)\leq c(X,Y_\ell)+b_\ell$ for all $X\in \Omega$ and $\phi(Z_{m_j})= c(Z_{m_j},Y_\ell)+b_\ell$.
Since $Z_m\to Z$ as $m\to \infty$, it follows by continuity that $\phi(Z)= c(Z,Y_\ell)+b_\ell$, and so $Y_\ell \in \mathcal N_\phi(Z)$.
On the other hand, since $Z\in \mathcal T_\phi(F)$ we obtain $Z\in \mathcal Q$.
\end{proof}

We are now in a position to prove the following comparison principle, akin to \cite[Theorem 2.7]{gutierrez-huang:nearfieldrefractor}, which clearly implies the uniqueness in
Theorem \ref{thm:existence of metasurfaces near field}.

\begin{proposition}\label{prop:comparison principle}
Let $\b=(b_1,\cdots ,b_N)$, $\b^*=(b_1^*,\cdots ,b_N^*)$, and let $\phi_{\b}, \phi_{\b^*}$ be given by \eqref{eq:solutionsecondbdryMA near field} two admissible phases solving \eqref{eq:energy conservation condition local} with the density $\rho>0$ a.e; and assume $\Omega$ is connected. 

If $b_1^*\geq b_1$, then $b_i^*\geq b_i$ for all $1\leq i\leq N$.
So if $b_1^*= b_1$, then $b_i^*= b_i$ for all $1\leq i\leq N$.
\end{proposition}
\begin{proof}
Let $I=\{i:b_i^*\geq b_i\}$ and $J=\{j:b_j> b_j^*\}$. We then want to show that $J=\emptyset$.
Suppose by contradiction that $J\neq \emptyset$.
From \eqref{eq:energy conservation condition local}, $\int_{\mathcal T_{\phi_\b}(Y_i)}\rho(X)\,dX=g_i>0$ for $1\leq i\leq N$. Hence $\mathcal T_{\phi_\b}(Y_i)$ has positive measure for each $1\leq i\leq N$,
and likewise $\mathcal T_{\phi_{\b^*}}(Y_i)$.

Set $F=\{Y_j:j\in J\}$. We shall prove that
\begin{equation}\label{eq:inclusion of closure into interior}
\overline{\mathcal T_{\phi_\b}\(F\)^\circ}
\subset 
\mathcal T_{\phi_{\b^*}}\(F\)^\circ.
\end{equation}

Let $X\in \overline{\mathcal T_{\phi_\b}\(F\)^\circ}$. 
We first claim that there exists $Y_j\in F$ such that $c(\cdot, Y_j)+b_j$ supports $\phi_\b$ at $X$.
Indeed, for each $m\geq 1$, $B_{1/m}(X)\cap \mathcal T_{\phi_\b}\(F\)^\circ\neq \emptyset$, where $B_{1/m}(X)$ is the open ball with radius $1/m$ centered at $X$. 
Since the last intersection is a non empty open set, it has positive measure. Hence $|B_{1/m}(X)\cap \mathcal T_{\phi_\b}\(F\)^\circ\cap \mathcal Q^c|>0$ for all $m$ where 
$\mathcal Q$ is the null set where $\mathcal N_{\phi_\b}$ is not a singleton.
Let us then pick $Z_m\in B_{1/m}(X)\cap \mathcal T_{\phi_\b}\(F\)^\circ\cap \mathcal Q^c$ and proceed as in the proof of Lemma \ref{lm:boundary of tracing contained in singular set}.
That is, taking a subsequence $\mathcal N_{\phi_\b}(Z_{m_\ell})$ equals one value $Y_j$ in $F$, so $c(\cdot,Y_j)+b_j$ supports $\phi_\b$ at $Z_{m_\ell}$. Letting $\ell\to \infty$ we obtain that $c(\cdot,Y_j)+b_j$ supports $\phi_\b$ at $X$.
The claim is then proved.
 

To prove \eqref{eq:inclusion of closure into interior}, we then have 
\begin{align*}
c(X,Y_i)+b_i^*&\geq c(X,Y_i)+b_i\quad \forall i\in I\\
&\geq \phi_{\b}(X)=c(X,Y_j)+b_j,\quad \text{since $X\in \overline{\mathcal T_{\phi_\b}\(F\)^\circ}$}\\
&>c(X,Y_j)+b_j^*\quad \text{since $j\in J$,}
\end{align*} 
so $\min_{i\in I}c(X,Y_i)+b_i^*>c(X,Y_j)+b_j^*$. 
Then by continuity of $c$, there exists an open neighborhood $N_X$ of the point $X$ such that 
\[
\min_{i\in I}c(Z,Y_i)+b_i^*>c(Z,Y_j)+b_j^*\quad \forall Z\in N_X.
\]
Hence 
\[
\phi_{b^*}(Z)=\min_{1\leq i\leq N}c(Z,Y_i)+b_i^*=\min_{i\in J}c(Z,Y_i)+b_i^*\quad \text{for all $Z\in N_X$.}
\]
This implies that given $Z\in N_X$ there exists $m\in J$, depending on $Z$, such that 
$\phi_{b^*}(Z)=c(Z,Y_m)+b_m^*$. By definition $\phi_{b^*}(Y)\leq c(Y,Y_m)+b_m^*$ for all $Y\in \Omega$.
Then $c(\cdot,Y_m)+b_m^*$ supports $\phi_{b^*}$ at $Y=Z$, that is, $Z\in \mathcal T_{\phi_{\b^*}}(Y_m)$.

Therefore for each $X\in \overline{\mathcal T_{\phi_\b}\(F\)^\circ}$ there exists an open neighborhood $N_X$ of $X$ such that 
\[
N_X\subset \cup_{m\in J}\mathcal T_{\phi_{\b^*}}(Y_m)=\mathcal T_{\phi_{\b^*}}\(F\)
\]
which proves \eqref{eq:inclusion of closure into interior}.

Since $\mathcal T_{\phi_{\b^*}}\(F\)$ is closed, $\mathcal T_{\phi_{\b^*}}\(F\)=\mathcal T_{\phi_{\b^*}}\(F\)^\circ \cup \partial \(\mathcal T_{\phi_{\b^*}}\(F\)\)$. Then from Lemma \ref{lm:boundary of tracing contained in singular set} and since $\rho>0$ a.e., we obtain 
\[
0<\sum_{j\in J}g_j=\int_{\mathcal T_{\phi_{\b^*}}(F)}\rho(x)\,dx
=
\int_{\mathcal T_{\phi_{\b^*}}(F)^\circ}\rho(x)\,dx,
\]
in particular, the set $\mathcal T_{\phi_{\b^*}}(F)^\circ\neq \emptyset$.
From \eqref{eq:inclusion of closure into interior} and connectedness of $\Omega$, then set $\mathcal T_{\phi_{\b^*}}\(F\)^\circ
\setminus
\overline{\mathcal T_{\phi_\b}\(F\)^\circ} 
$
is a non empty open set and therefore it has positive measure.
Obviously, $\mathcal T_{\phi_\b}\(F\)^\circ\subset \overline{\mathcal T_{\phi_\b}\(F\)^\circ}$, and so
$\mathcal T_{\phi_{\b^*}}\(F\)^\circ
\setminus
\mathcal T_{\phi_\b}\(F\)^\circ
$
also has positive measure. 
Since $\rho>0$ a.e., we obtain
\[
\int_{\mathcal T_{\phi_{\b^*}}(F)^\circ}\rho(x)\,dx>\int_{\mathcal T_{\phi_{\b}}(F)^\circ}\rho(x)\,dx
\]
which is a contradiction since both sides of this inequality equal $\sum_{j\in J}g_j$.

Therefore $J=\emptyset$ which completes the proof of the proposition.
\end{proof}

\begin{remark}\rm
In Definition \ref{def:admissible phase near field}, the admissible phase $\phi$ is supported by the functions $c(\cdot,\cdot)+b$ from above which yields the concave-infinity case used to prove Theorem \ref{thm:existence of metasurfaces near field}. 
An alternative definition of admissible phase can be made with supporting functions $c(\cdot,\cdot)+b$ from below which yields the convex-infinity case. With this definition, proceeding in a similar way and using the results from
\cite[Sect. 2]{gutierrez-huang:nearfieldrefractor}, a theorem similar to \ref{thm:existence of metasurfaces near field} also follows, where in \eqref{eq:solutionsecondbdryMA near field} the min is replaced by max.
A reason to choose the Definition \ref{def:admissible phase near field}, is that this notion is more suitable for the initialization of the numerical scheme developed in Section \ref{sec:computation of the Laguerre cells} to compute the Laguerre cells.
\end{remark}

%% file: section3-mapping.tex
\section{Analysis of the refracted distribution}
\label{sec:kantorovitch formulation}

\begin{definition}[Laguerre cells and refracted distribution] We define the
Laguerre cells associated to $(Y_i,b_i)_{1\leq i\leq N}$ by
\begin{equation}\label{eq:definition of Laguerre cell}
\Lag_i(\b)=\left\{X\in \Omega: c(X,Y_i)+b_i\leq c(X,Y_k)+b_k\text{ for all $1\leq k\leq N$} \right\},\quad 1\leq i\leq N,
\end{equation} 
where $c$ is given in \eqref{eq:definition of the cost cXY}.  We call
\emph{refracted distribution} to the vector $G(\psi)\in\R^N$ defined by
\begin{equation}\label{eq:refractor mapping}
G(\b)=\(G_1(\b),\cdots ,G_N(\b)\), \quad  \hbox{ where } G_i(\b)=\int_{\Lag_i(\b)}\rho(X)\,dX,
\end{equation}
which encodes the amount of light refracted in each direction
$\{Y_1,\hdots,Y_N\}$.
\end{definition}

Given $\b=\(b_1,\cdots ,b_N\)$, let $\phi_{\b}$ be the function defined by \eqref{eq:solutionsecondbdryMA near field}, and let $\mathcal M_{\phi_\b}$ be the refracted measure defined by the left hand side of \eqref{eq:energy conservation condition local}, that is,
\[
\mathcal M_{\phi_\b}(E)=
\int_{\mathcal T_{\phi_\b}(E)}\rho(X)\,dX.
\]
One can easily verify that $\Lag_i(\b)=\mathcal T_{\phi_\b}(Y_i)$, so
that the refractor measure is given by
\[
\mathcal M_{\phi_\b} = \sum_{1\leq i\leq N} G_i(\b) \delta_{Y_i}.
\]
Assuming that $\mu = \sum_{1 \leq i\leq N} g_i \delta_{Y_i}$ with
$\mathbf{g}\in\R^N$, the near-field metasurface refractor problem
\eqref{eq:energy conservation condition local} means to solve the finite-dimensional non-linear
system of equations
\begin{equation} \label{eq:refractor-eqn}
  G(\b) = \g,
\end{equation}
where $\g = (g_1,\hdots,g_N)$. 
The goal of this section is to gather a few
properties of the refracted distribution, which will be used to establish
the global convergence of a damped Newton algorithm to solve this
system of equations.

\subsection{Regularity of the map $G$}

\begin{theorem}[Partial derivatives of $G$] \label{thm:pd}
  Let $0<\alpha <\beta$, let $\Omega \subseteq \R^2\times \{\alpha\}$
  be a polygon and let $\rho\in C^0(\Omega)$. Assume that the target
  $T=\{Y_1,\cdots ,Y_N\}$ is included in the plane $\{x_3=\beta\}$.
  Then the refracted distribution map $G$ given by \eqref{eq:refractor mapping}
  belongs to $C^1(\R^N)$ and the partial derivatives of $G$ are given
  by
\begin{equation}\label{eq:derivatives of Gi with respect bj j not equal i}
\dfrac{\partial G_i}{\partial b_j}(\b)
= G_{ij}(\b):=
\int_{\Lag_{ij}(\b)} \dfrac{\rho(X)}{\left| \nabla_Xc(X,Y_i)-\nabla_Xc(X,Y_j)\right|}\,dX\quad i\neq j,
\end{equation}
when $j\neq i$.
where the integration is over the curve
$$ \Lag_{ij}(\b)=\Lag_i(\b)\cap \Lag_j(\b). $$
The diagonal partial derivatives are given by
\begin{equation}\label{eq:partial derivative of G with respect to i}
\dfrac{\partial G_i}{\partial b_i}(\b)
= G_{ii}(\b) := -\sum_{j\neq i }\dfrac{\partial G_i}{\partial b_j}(\b),
\end{equation}
\end{theorem}

\begin{remark}[Assumptions on $\Omega$]\label{rmk:assumptions on Omega}
From the proof of the theorem, one can verify that the hypothesis that
$\Omega$ is a polygon can be replaced by the following two hypothesis:
  \begin{itemize}
    \item the boundary of $\Omega$ in $\R^2\times \{\alpha\}$ has area zero;
    \item the intersection of $\partial \Omega$ with any conic in $\R^2\times \{\alpha\}$ is finite.
  \end{itemize}
  Moreover, if $\rho$ has compact support in $\Omega$, then the
  assumption that $\partial \Omega$ has measure zero is not needed.
\end{remark}

Theorem~\ref{thm:pd} is proved in the same way as 
\cite[Thm.~45]{2019-merigot-thibert-optimal-transport}, provided that we are
able to show that the functions $G_{ij}$ are continuous, which 
replaces \cite[Lemma 46]{2019-merigot-thibert-optimal-transport}. Note
that unlike in Lemma 46, we do not need the points $Y_1,\hdots,Y_N$ to be in a generic position (\cite[Definition 16]{2019-merigot-thibert-optimal-transport}). 

\begin{lemma} Under the assumptions of Theorem~\ref{thm:pd}, the functions $G_{ij}$ are continuous.
\end{lemma}

\begin{proof} Define $H_{ij}(\b)=\{X\in \Omega:c(X,Y_i)+b_i=c(X,Y_j)+b_j\}$,
so that 
\begin{equation}\label{eq:inclusion of Laguerre cell in H}
\Lag_{ij}(\b)\subset H_{ij}(\b),
\end{equation}
$1\leq i,j\leq N$. Set
$f(X)=c(X,Y_i)-c(X,Y_j)$ for $X\in \Omega$, and let $\Omega\Subset
\Omega_1\Subset \Omega_2$ with $\Omega_i$ open sets and let $\bar f$
be an extension of $f$ to $\Omega_2$ so that the gradient of $\bar f$
is continuous in $\Omega_2$ and agrees with the gradient of $f$ in
$\Omega$.  Consider the system of ODEs
\[
\begin{cases}
&\dfrac{d\Phi}{dt}(t,X)=\dfrac{\nabla \bar f\(\Phi(t,X)\)}{\left|\nabla \bar f\(\Phi(t,X)\)\right|^2}:=F\(\Phi(t,X)\)\\
&\Phi(0,X)=X.
\end{cases}
\]
For $\epsilon$ sufficiently small, there is a unique local solution $\Phi:(-\epsilon,\epsilon)\times \Omega_1\to \Omega_2$ and so that $\Phi\((-\epsilon,\epsilon)\times \Omega\)\subset \Omega_1$. Here $F(X)$ is $C^1$ for $X$ in $\Omega_2$.
Since $\dfrac{d}{dt}\bar f\(\Phi(t,X)\)=1$, we have for $-\epsilon<t<\epsilon$ and $X\in \Omega_1$
\begin{equation}\label{eq:formula for f}
\bar f\(\Phi(t,X)\)=\bar f\(\Phi(0,X)\)+t.
\end{equation}
Since $F$ is $C^1(\Omega_2)$, we also have that $\Phi(t,X)\to X$ in $C^1(\Omega_1)$ as $t\to 0$.

Now, let $\b^n=\(b_1^n,\cdots ,b_N^n\)\to \b=\(b_1,\cdots ,b_N\)$ as $n\to \infty$.
We shall prove that $G_{ij}(\b^n)\to G_{ij}(\b)$ as $n\to \infty$.
Set 
\[
a:=b_j-b_i,\quad t_n:=b_j^n-b_i^n-a,\quad L:=\Lag_{ij}(\b),\quad L_n:=\Phi\(-t_n,\Lag_{ij}\(\b^n\)\).
\]
Let $H=\bar f^{-1}(a)=\{X\in \Omega_1:\bar f(X)=a\}$. From \eqref{eq:inclusion of Laguerre cell in H}, $L\subset H\cap \Omega$. Also, if $X\in L_n$, there is $Y\in \Lag_{ij}\(\b^n\)\(\subset H_{ij}(\b^n)\)$ such that $X=\Phi(-t_n,Y)$,  and 
from \eqref{eq:formula for f} $\bar f(X)=\bar f(Y)-t_n=b_j^n-b_i^n-t_n=a$. Since $Y\in \Omega$ we have $X\in \Omega_1$. Thus, $L_n\subset H$ for $n$ large since $t_n\to 0$.
Define for $X\in H$,
\[
F_n(X)=\Phi(-t_n,X),
\]
$F_n: H\to \Phi(-t_n,H)\subset \Omega_2$.
We have $L_n=F_n\(\Lag_{ij}\(\b^n\)\)\subset \Omega_1$.
Write
\begin{align*}
G_{ij}(\b^n)&=\int_{\Lag_{ij}(\b^n)} \dfrac{\rho(Y)}{\left| \nabla_Xc(Y,Y_i)-\nabla_Xc(Y,Y_j)\right|}\,dY\\
&=\int_{F_n^{-1}(L_n)} \dfrac{\rho(Y)}{\left| \nabla_Xc(Y,Y_i)-\nabla_Xc(Y,Y_j)\right|}\,dY.
\end{align*}
Since $\rho$ is defined in $\Omega$ and in the last integral we will make the change of variables $Y=F_n(X)$\footnote{We notice that $F_n$ is a genuine change of variables because letting $\dfrac{\partial \Phi}{\partial X}(t,X)$ being the Jacobian matrix of $\Phi$, setting $J(t,X)=\det \(\dfrac{\partial \Phi}{\partial X}(t,X)\)$ we have that $J$ satisfies the following ode: 
$$\dfrac{\partial J}{\partial t}(t,X)=\diver(F)\(\Phi(t,X)\)\,J(t,X)$$ 
with the initial condition $J(0,X)=1$.}, we need to extend $\rho$ to $\Omega_2$, so let $\bar \rho$ be a bounded extension of $\rho$ to $\Omega_2$ so that $\bar \rho\in C(\Omega_2)$.
So
\begin{align*}
G_{ij}(\b^n)&=\int_{L_n} \dfrac{\rho(F_n(X))}{\left| \nabla_Xc(F_n(X),Y_i)-\nabla_Xc(F_n(X),Y_j)\right|}\,|J_{F_n}(X)|\,\,dX\\
&=\int_{H} \dfrac{\bar \rho(F_n(X))}{\left| \nabla_Xc(F_n(X),Y_i)-\nabla_Xc(F_n(X),Y_j)\right|}\,|J_{F_n}(X)|\,\chi_{L_n}(X)\,dX.
\end{align*}
Now $F_n(X)\to \Phi(0,X)=X$ and $|J_{F_n}(X)|\to 1$ as $n\to \infty$.
Hence to show that $G_{ij}(\b^n)\to G_{ij}(\b)$, is enough to show that $\chi_{L_n}(X)\to \chi_L(X)$ for a.e. $X$.
Given an arbitrary sequence of sets $E_n$ we have the following
\begin{align*}
\limsup_{k\to \infty}E_k&=\cap_{n=1}^\infty\cup_{k\geq n}E_k\\
\liminf_{k\to \infty}E_k&=\cup_{n=1}^\infty\cap_{k\geq n}E_k\\
\chi_{\limsup_{k\to \infty}E_k}(X)&=\limsup_{n\to \infty}\chi_{E_n}(X)\\
\chi_{\liminf_{k\to \infty}E_k}(X)&=\liminf_{n\to \infty}\chi_{E_n}(X).
\end{align*}
We first prove that $\chi_{\limsup_{n\to \infty}L_n}(X)\leq \chi_{L}(X)$,
which is equivalent to show that $\limsup_{n\to \infty}L_n\subset L$.
In fact, if $X\in \limsup_{n\to \infty}L_n$, there exists a subsequence $n_\ell$ such that $X\in L_{n_\ell}$ for $\ell=1,2,\cdots$. So $X=\Phi(-t_{n_\ell},Z_{n_\ell})$ for some $Z_{n_\ell}\in \Lag_{ij}(\b^{n_\ell})$. Since $\Omega$ is compact, there is a subsequence $Z_{n_{\ell_m}}\to Z\in \Lag_{ij}(\b)=L$. Hence $X=\Phi(0,Z)=Z$, so $X\in L$ as desired.

Let 
\begin{equation}\label{eq:singular set}
S=\bigcup_{k\neq i,j}H_{ijk}(\b)\cup \(\Lag_{ij}(\b)\cap \partial \Omega\).
\end{equation}
Under the assumptions on the boundary $\partial \Omega$ described below, we shall prove in a moment that this is a set of linear measure zero.
Taking this for granted, we claim that 
\begin{equation}\label{eq:inclusion L into lim inf}
L\subset \(\liminf_{n\to \infty}L_n\)\cup S
\end{equation}
obtaining $\chi_L(X)\leq \liminf_{n\to \infty}\chi_{L_n}(X)$ for $X\notin S$. 
Therefore the sequence $\chi_{L_n}(X)\to \chi_L(X)$ for a.e. $X$. 
Since $\bar \rho$ is continuous and bounded, we then obtain by Lebesgue dominated convergence theorem that
\begin{align*}
G_{ij}(\b^n)&=\int_{H} \dfrac{\bar \rho(F_n(X))}{\left| \nabla_Xc(F_n(X),Y_i)-\nabla_Xc(F_n(X),Y_j)\right|}\,|J_{F_n}(X)|\,\chi_{L_n}(X)\,dX\\
&\qquad \to 
\int_{H} \dfrac{\bar \rho(X)}{\left| \nabla_Xc(X,Y_i)-\nabla_Xc(X,Y_j)\right|}\,\chi_{L}(X)\,dX=G_{ij}(\b)
\end{align*}
as $n\to \infty$ showing that $G_{ij}$ is continuous.

It then remains to prove the claim \eqref{eq:inclusion L into lim inf}. 
Indeed, if $X\in L$ and $X\in \partial \Omega$, then $X\in S$.
On the other hand, if $X\in L\cap \text{interior}(\Omega)$, and $X\notin S$, then we show $X\in \liminf_{n\to \infty}L_n$. We have
\[
X\in S^c=\bigcap_{k\neq i,j}H_{ijk}(\b)^c\cap \(\Lag_{ij}(\b)\cap \partial \Omega\)^c.
\] 
Since $H_{ijk}(\b)=H_{ij}(\b)\cap H_{ik}(\b)$, we then have 
\[
X\in \bigcap_{k\neq i,j}\(H_{ij}(\b)^c\cup H_{ik}(\b)^c\)=H_{ij}(\b)^c \cup \bigcap_{k\neq i,j}H_{ik}(\b)^c.
\]
Since $X\in L=\Lag_{ij}(\b)\subset H_{ij}(\b)$, we get 
\[
X\in \bigcap_{k\neq i,j}H_{ik}(\b)^c,
\]
that is, $c(X,Y_i)-c(X,Y_k)\neq b_k-b_i$ for all $k\neq i$. On the other hand, since $X\in \Lag_i(\b)$ we have $c(X,Y_i)-c(X,Y_k)<b_k-b_i$ for all $k\neq i$.
Since $\Phi(t_n,X)\to X$ and $\b^n\to \b$, it follows that there exists $n_0$ such that 
\[
c\(\Phi(t_n,X),Y_i\)-c\(\Phi(t_n,X),Y_k\)<b_k^n-b_i^n\quad \forall k\neq i
\]
for all $n\geq n_0$. That is, $\Phi(t_n,X)\in \Lag_{ij}(\b^n)$ for all $n\geq n_0$.
Now $X\in L_n$ iff $X=\Phi(-t_n,Y)$ for some $Y\in \Lag_{ij}(\b^n)$. But
$\Phi(t_n,X)=\Phi\(t_n, \Phi(-t_n,Y)\)=\Phi(t_n+(-t_n),Y)=\Phi(0,Y)=Y$ from the semigroup property of the flow. Therefore $X\in L_n$ iff $\Phi(t_n,X)\in \Lag_{ij}(\b^n)$, and the claim is then proved.

To complete the analysis we show that the set $S$ in \eqref{eq:singular set} has measure zero.
Indeed, recall that the target $T$ is contained in the plane $x_3=\beta$ and $\Omega$ is contained in the plane $x_3=\alpha$ with $\alpha<\beta$. We claim that the set of points $X$ with $x_3=\alpha$ satisfying 
\begin{equation}\label{eq:degeneracy condition}
\begin{cases}
&c(X,Y_0)-c(X,Y_1)=t_1\\
&c(X,Y_0)-c(X,Y_2)=t_2
\end{cases}
\end{equation}
is a discrete set for any $Y_0,Y_1,Y_2$ distinct points in $x_3=\beta$ and for all $t_1,t_2\in \R$.
In our case \eqref{eq:degeneracy condition} reads
\[
\begin{cases}
&|X-Y_0|-|X-Y_1|=t_1\\
&|X-Y_0|-|X-Y_2|=t_2
\end{cases}
\]
Each of these equations describe a hyperboloid of two sheets, one with foci $Y_0,Y_1$, and the other with foci $Y_0,Y_2$. These two hyperboloids are intersected with the plane $x_3=\alpha$ where $X$ lies.
Since hyperboloids are quadric surfaces, their intersection with the plane $x_3=\alpha$ are conics. 
Now, two conics in the plane intersect in a finite number of points unless they are equal. But if they are equal, then the foci $Y_1$ and $Y_2$ must be the equal which is impossible.
This shows that the first set in the union \eqref{eq:singular set} has measure zero.
Finally we note that from \eqref{eq:inclusion of Laguerre cell in H} and using that $\partial \Omega$ is a polygon,  the intersection $\partial \Omega\cap H_{ij}(\b)$ is finite. This implies that the second set $\Lag_{ij}(\b)\cap \partial \Omega$ also has measure zero.
\end{proof}


\subsection{Monotonicity of the map $G$}

Denote $DG(\b)$ the Jacobian matrix of $G$ at a point $\b\in\R^N$.
By invariance of the Laguerre cells under addition of a constant, the
one-dimensional space $\R \e$, with $\e=(1,1,\cdots ,1)$, is always
included in $\ker(DG(\b))$.  The next theorem proves the converse
inclusion, i.e. $\ker(DG(\b)) = \R e$, whenever all the Laguerre
cells have positive mass, i.e. $G_i(\b) >0$ for all
$i\in\{1,\hdots,N\}$. This implies a strong monotonicity of the
refracted distribution map $G$, which is used to prove convergence of a
damped Newton algorithm in the next section.

\begin{theorem} \label{thm:monotonicity}
  Assume that the conditions of Theorem~\ref{thm:pd} hold, that $\rho\geq 0$, and that 
  $$\mathbf{Z}=\text{\rm Int}(\Omega)\cap\{X\in \Omega:\rho(X)>0\} $$
  is a connected set. Let $G$ be the mapping given in \eqref{eq:refractor mapping},
let 
\[
S_+=\{\b=(b_1,\cdots ,b_N): G_i(\b)>0\quad \forall 1\leq i\leq N\},
\]
and set $\e=(1,1,\cdots ,1)$.  Then for each $\b\in S_+$ the matrix
$DG(\b)$ is symmetric non-positive definite and
$\text{ker}\(DG(\b)\)=\R\,\e$.
\end{theorem}

\begin{definition}
The $N\times N$ matrix $H$ is reducible if there exist non empty sets $A,B$ with $A\cap B=\emptyset$, $A\cup B=\{1,\cdots ,N\}$ such that $H_{\alpha\beta}=0$ for $\alpha\in A$ and $\beta\in B$.
The matrix $H$ is irreducible if it is not reducible.
\end{definition}

\begin{lemma}\label{lm:irreducible matrices}
Let $H = (H_{ij})_{1\leq i\leq j}$ be an $N\times N$ symmetric matrix satisfying $H_{i,j} \geq 0$ for $i\neq j$ and 
\begin{equation}\label{eq:sum of rows equals zero}
H_{ii}=-\sum_{j\neq i}H_{ij}\quad 1\leq i\leq N.
\end{equation}
Then $H$ is negative semidefinite.
If in addition $H$ is irreducible, then $\text{ker}\(H\)=\R e$.
\end{lemma}




\begin{proof}[Proof of Theorem~\ref{thm:monotonicity}]
Here we apply the previous Lemma to $H=DG(\b)$. 
From \eqref{eq:partial derivative of G
  with respect to i}, \eqref{eq:sum of rows equals zero} holds, and
from \eqref{eq:derivatives of Gi with respect bj j not equal i},
$G_{ij}(\b)\geq 0$ for any $j\neq i$. To prove the theorem, by Lemma
\ref{lm:irreducible matrices} it suffices to prove  that the matrix
$DG(\b):=\(G_{ij}(\b)\)_{i,j}$ is irreducible for $\b\in S_+$.  We proceed in steps.

{\bf Step 1.}
If 
\[
S=\bigcup_{Y_i\neq Y_j\neq Y_k}\Lag_{ijk}(\b),
\]
then $\mathbf{Z}\setminus S$ is open and path-connected.  Indeed, by
the proof of Theorem~\ref{thm:pd}, we know that the set $S$ has zero
length, or more precisely zero one-dimensional Hausdorff measure. By
\cite[Lemma 49]{2019-merigot-thibert-optimal-transport}, the fact that
$\mathbf{Z}$ is path-connected implies that $\mathbf{Z}\setminus S$ is
path-connected.

{\bf Step 2.}
For each $1\leq i\leq N$, $\Lag_i(\b)\cap \(\mathbf{Z}\setminus S\)\neq \emptyset$.

Since $\b\in S_+$, $\rho\(\Lag_i(\b)\) := \int_{\Lag_i(\b)} \rho(X)\, dX >0$.
Since $S$ is negligible,
$$ \rho\(\Lag_i(\b)\cap \(\mathbf{Z}\setminus S\)\)
=\rho\(\Lag_i(\b)\cap \mathbf{Z}\) = \rho(\Lag_i(\b)),$$
where we used the definition of $\mathbf{Z}$ and the assumption $\rho(\partial\Omega) = 0$ to get the last equality.  Since by assumption $\rho(\Lag_i(\b)) > 0$, we directly get that $\Lag_i(\b)\cap \(\mathbf{Z}\setminus S\)$ is nonempty.

%

{\bf Step 3.}
Let $i\neq j$. If $X\in \(\mathbf{Z}\setminus S\)\cap \Lag_i(\b)\cap \Lag_j(\b)$, then $G_{ij}(\b)>0$.

Indeed, since $X\in S^c$, 
\begin{align*}
c(X,Y_i)+b_i&=c(X,Y_j)+b_j\\
c(X,Y_i)+b_i&<c(X,Y_k)+b_k\quad \forall k\neq i,j.
\end{align*}
Then by continuity of $c$ and $\rho$ there exists a ball $B_r(X)$ such that 
\[
c(X',Y_i)+b_i<c(X',Y_k)+b_k\quad \forall X'\in B_r(X)\quad \forall k\neq i,j,
\]
and with
\[\rho(X')>0\quad \forall X'\in B_r(X).
\]
This implies that
\[
\left\{X':c(X',Y_i)+b_i=c(X',Y_j)+b_j\right\}\cap B_r(X)\subset \Lag_i(\b)\cap \Lag_j(\b) =\Lag_{ij}(\b).
\]
As shown in the proof of Theorem~\ref{thm:pd}, the set $\left\{X':c(X',Y_i)+b_i=c(X',Y_j)+b_j\right\}$ is a conic and in particular a $1$-dimensional manifold.
Therefore
\begin{align*}
G_{ij}(\b)&=\int_{\Lag_{ij}(\b)} \dfrac{\rho(X')}{\left| \nabla_Xc(X',Y_i)-\nabla_Xc(X',Y_j)\right|}\,dX'\\
&\geq 
\int_{\left\{X':c(X',Y_i)+b_i=c(X',Y_j)+b_j\right\}\cap B_r(X)} \dfrac{\rho(X')}{\left| \nabla_Xc(X',Y_i)-\nabla_Xc(X',Y_j)\right|}\,dX'>0.
\end{align*}

To conclude the proof of the irreducibility of $DG$, we  suppose by contradiction that there exists $\b\in S_+$ such that the matrix $DG(\b)$ is reducible. This means there exist non empty disjoint sets $I$ and $J$ such that $\{1,\cdots ,N\}=I\cup J$ with $G_{ij}(\b)=0$ for all $(i,j)\in I\times J$.
Let 
\[
Z_I=\bigcup_{i\in I}\Lag_i(\b)\cap (\mathbf{Z}\setminus S),\quad Z_J=\bigcup_{i\in J}\Lag_i(\b)\cap (\mathbf{Z}\setminus S).
\]
Then from Steps 2 and 3, the sets $Z_I$ and $Z_J$ are non empty and disjoint. Since $\cup_{1\leq i\leq N}\Lag_i(\b)=\Omega$ and $\Lag_i(\b)$ are closed in $\Omega$, we have that $\Lag_i(\b)\cap (\mathbf{Z}\setminus S)$ are relatively closed subsets of $\mathbf{Z}\setminus S$ with $\mathbf{Z}\setminus S=Z_I\cup Z_J$ contradicting the connectedness of $\mathbf{Z}\setminus S$.
\end{proof}



%% file: section4-numerics.tex
\section{Implementation and numerical experiments}\label{sec:computation of the Laguerre cells}
\subsection{Damped Newton algorithm}\label{sec:algorithm}
As shown in the previous section, solving the near-field metasurface refractor problem with target $\mu = \sum_{1\leq i\leq N} g_i \delta_{Y_i}$ amounts to solve the non-linear system $G(\b) = \g$ (see \eqref{eq:refractor-eqn}) with $\g = (g_1,\hdots,g_N)$. 
We use below the damped Newton algorithm introduced in~\cite{kitagawa-merigot-thibert} to solve this equation.
To do so, we pick an initialization vector $\b^0$ (see Remark \ref{rm:initialization vector}) such that all Laguerre cells have a positive amount of mass, and we denote
\[
\epsilon=
\dfrac12\,
\min\left\{ 
\min_{1\leq i\leq N}G_i(\b^0),\min_{1\leq i\leq N}g_i
\right\}>0
\]
Adding a constant to $\b^0$ if necessary, we may assume that $b^0 \in
\{e\}^\perp$ where $\e = (1,\hdots,1)$.  Thus, $\b^0$ belongs to the
set
$$ S = \{\b\in\R^N : \forall i\in\{1,\hdots,N\} G(\b) \geq \eps\} \cap \{\e\}^\perp. $$

\subsubsection{Algorithm}
Given an iterate $\b^k$, we explain how to define the next iterate
$\b^{k+1}$. We first denote $\mathbf{v}^k$ the solution to the system
\begin{equation}\label{eq:newt}
\begin{cases}
DG(\b)\,\mathbf{v}^k &= \g-G(\b)\\
\sum_{i=1}^N v^k_i&=1,
\end{cases}
\end{equation}
which exists and is unique by Theorem \ref{thm:monotonicity}.
Denoting $\b^{k}_\tau = \b^{k} + \tau \mathbf{v}^k$, we introduce
\begin{align*}
&E^k=\left\{2^{-\ell}:\ell\in \N, \b^k_{2^{-\ell}}\in S,
\left|G\(\b^k_{2^{-\ell}}\)-\g\right|\leq
\(1-\dfrac{2^{-\ell}}{2}\)\,|G(\b^k)-\g| \right\}\\
&\tau^k =  \max E^k.
\end{align*}
We then denote
$$ \b^{k+1} = b^k + \tau^k \mathbf{v}^k.$$

\begin{proposition}[Linear convergence]\label{prop:linearcv}
  Under the assumptions of Theorems~\ref{thm:pd} and
  \ref{thm:monotonicity}, there exists a constant $\tau^* \in (0,1)$
  such that
    $$ \nr{G(\b^k)-\g)} \leq \(1-\dfrac{\tau^*}{2}\)^k \nr{G(\b^0)-\g} $$
\end{proposition}

\begin{proof} 
Thanks to Theorems~\ref{thm:pd} and \ref{thm:monotonicity}, the
refracted distribution $G$ satisfies the assumptions of 
\cite[Prop.~50]{2019-merigot-thibert-optimal-transport}, which implies the result.
\end{proof}

\subsection{Computation of Laguerre cells}
\begin{figure}[H]
\begin{center}

\begin{tabular}{ccc}
\includegraphics[width=0.32\textwidth]{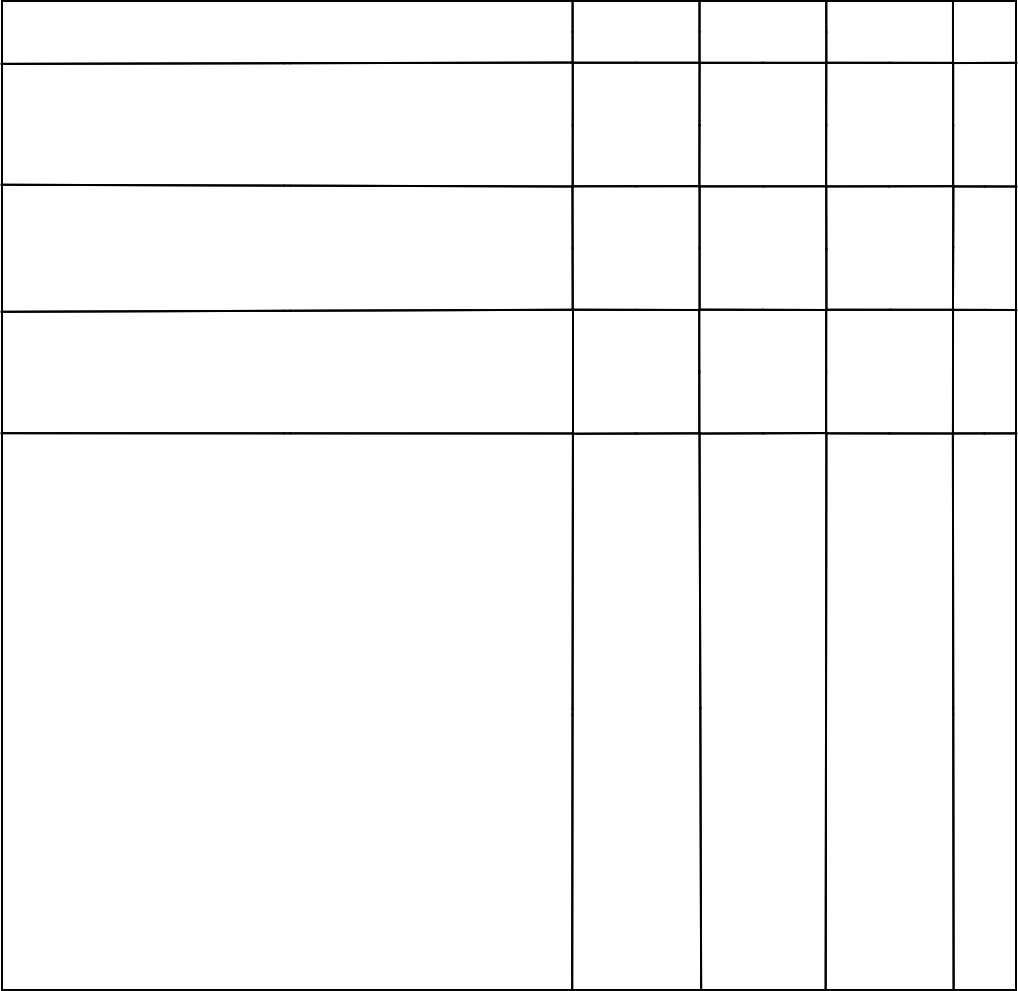}&
 \includegraphics[width=0.32\textwidth]{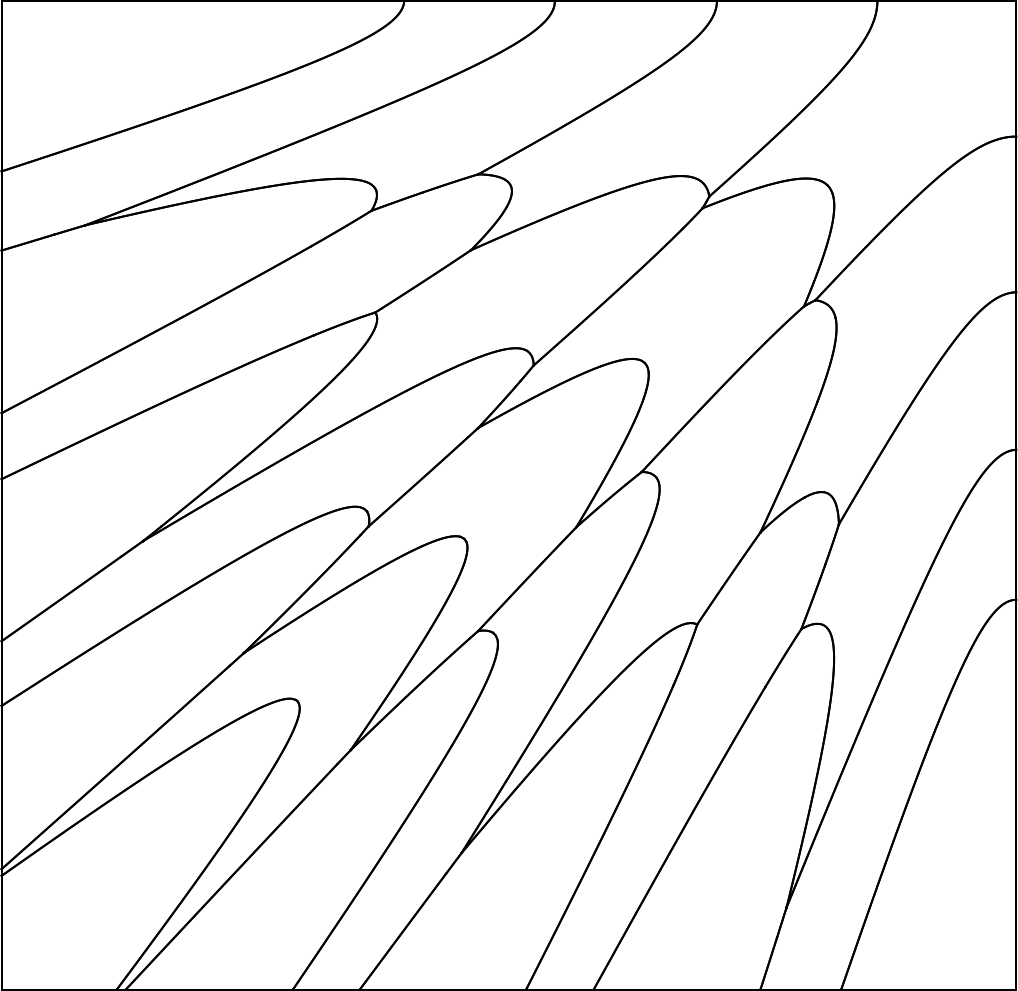}& \includegraphics[width=0.32\textwidth]{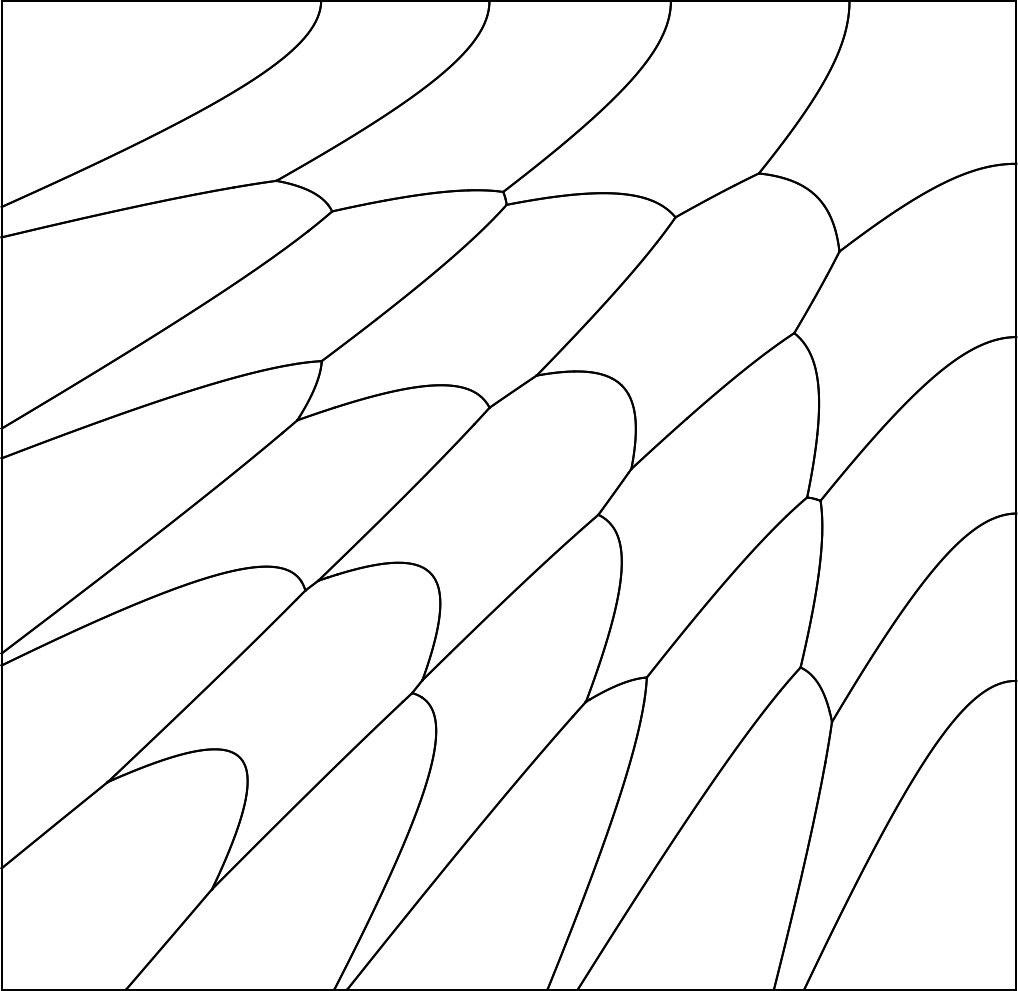}\\
\includegraphics[width=0.32\textwidth]{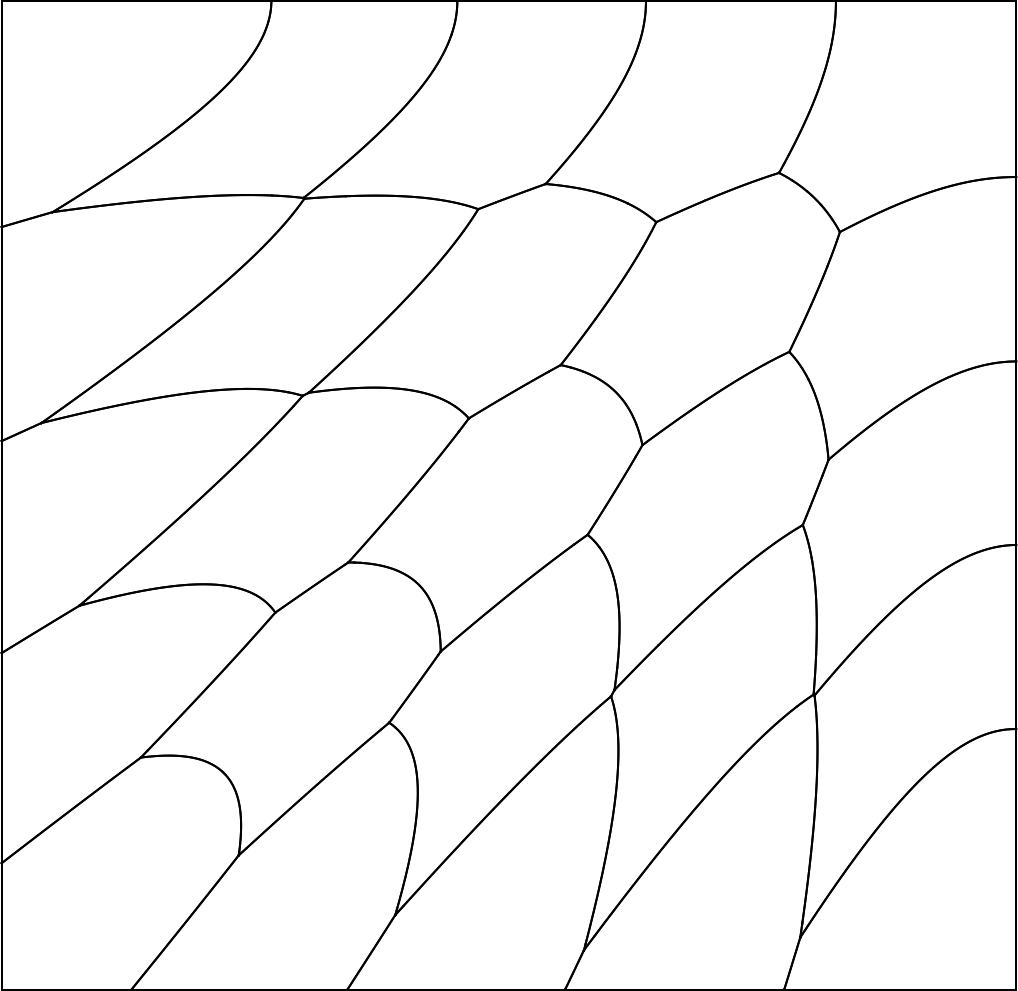}&
 \includegraphics[width=0.32\textwidth]{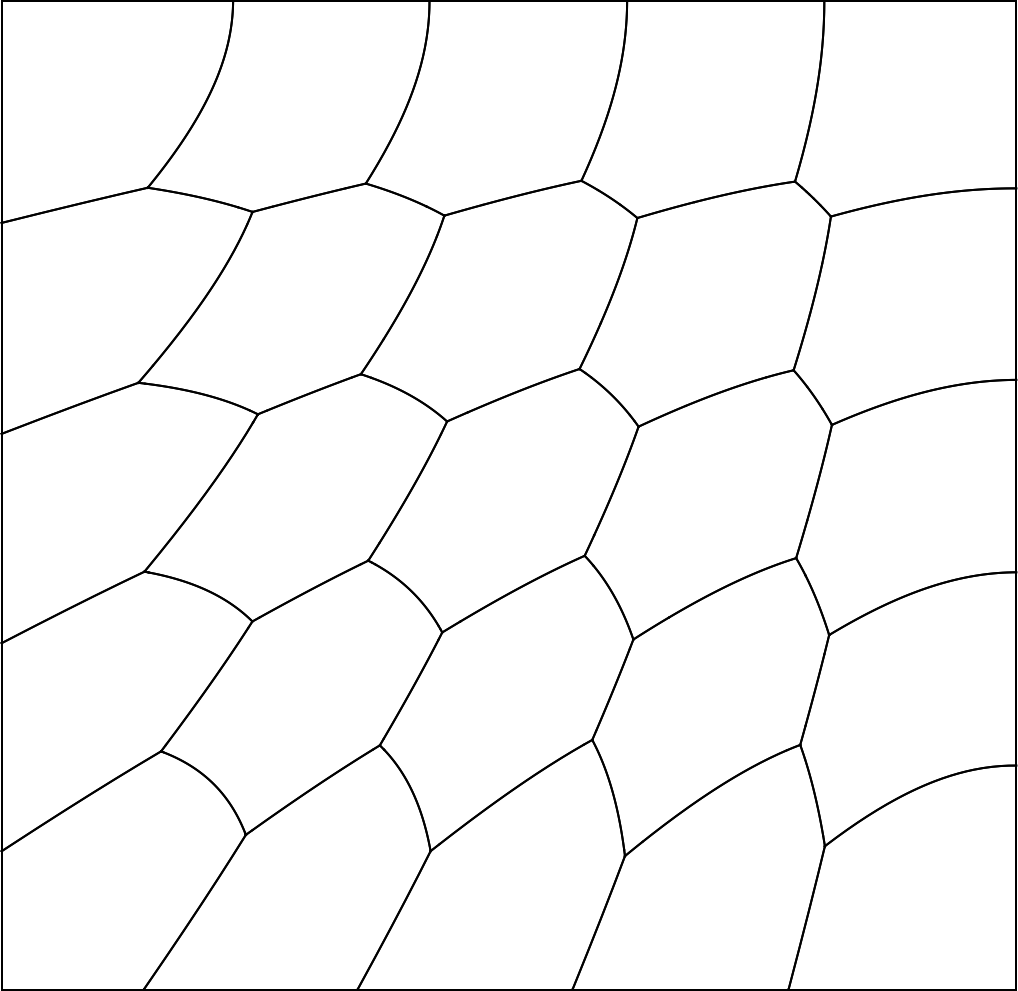}&
 \includegraphics[width=0.32\textwidth]{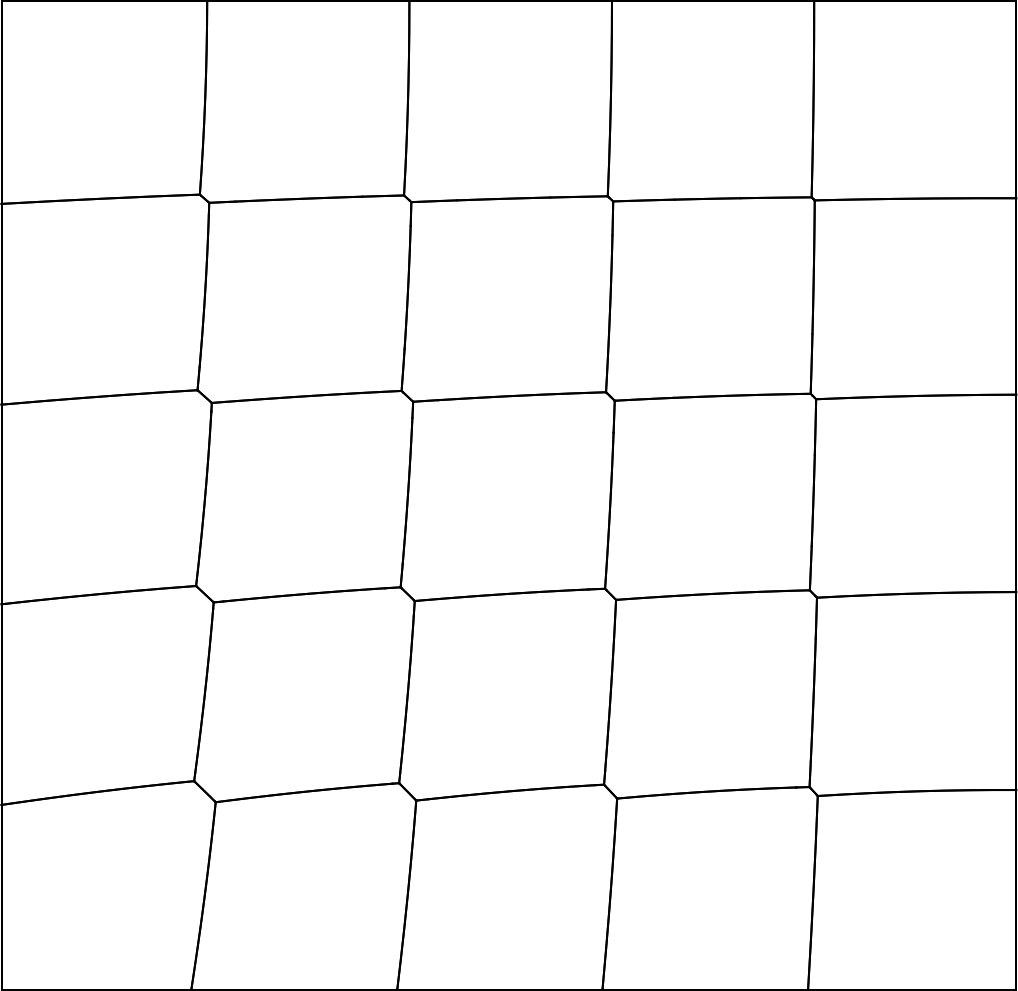}
 \end{tabular}
    \end{center}
    \caption{From left to right and top to bottom: The first image displays the Voronoi cells of $\{Y_1,\hdots,Y_N\}$. Then, each image displays the Laguerre cells associated to the solutions of the near-field metasurface problem described in \S\ref{subsec:numerics-beta} for $\delta \in \{0.1,0.2,0.3,0.5,2\}$. \label{fig:numerics-beta}}
\end{figure}

Computing the Laguerre cells associated to the near-field refractor
metasurface problem is not an easy task, because these cells have
curved boundaries (see Figure~\ref{fig:numerics-beta}), are not convex, etc.  In this section, we show that the Laguerre cells
can be obtained by using power diagrams. The advantage of this
formulation is that there are very efficient algorithms and
software libraries available to construct 3D power diagrams, with near-linear
complexity in $N$ for non-degenerate input.

\begin{definition}
Let $Q=\{(q_i,\omega_i)\}_{1\leq i \leq N}$ be a weighted cloud point set, i.e., $q_i \in \R^3$ and $\omega_i\in \R$. Then for each $1\leq i\leq N$, the $i$-th power diagram of $Q$ is defined by 
$$
\Pow_i(Q) = \left\{x\in \R^3: |x - q_i|^2 + \omega_i \leq |x - q_j|^2 + \omega_j \quad \forall j\in\{1,\cdots, N\} \right\}.
$$
\end{definition}
Let us define
$$
H_i(\b)=\left\{ X\in \R^2\times\{\alpha\}: |X-Y_i| + b_i \leq |X-Y_j|+b_j,\quad \forall 1\leq j\leq N\right\},\quad 1\leq i\leq N,
$$
with $\b=\(b_1,\cdots ,b_N\)$, $0 < \alpha<\beta$, and the points $Y_i=(y_i,\beta)$ lie in the horizontal plane $\R^2\times\{\beta\}$. 
From the definition of Laguerre cell in \eqref{eq:definition of Laguerre cell}, we have 
\[
\Lag_i(\b)=H_i(\b)\cap \Omega.
\]

\begin{proposition}[Point Source/Near-Field]\label{prop:computationlaguerre_PSNF}
We assume the following condition 
\begin{equation}\label{eq:bi and bj are not too far}
\forall i,j\in \{1,\cdots, N\} \quad |b_i - b_j| < \sqrt{4(\alpha-\beta)^2+|y_i-y_j|^2}.
\end{equation}
Then, for each $i \in \{1,\cdots, N\}$, one has
$$
H_i(\b)=\proj_{\R^2\times\{\alpha\}}\left( 
\Pow_i(Q) \cap \Sigma_i^+
\right)
$$
where   $\proj_{\R^2\times\{\alpha\}}$ denotes the orthogonal projection onto $\R^2\times\{\alpha\}$, $Q=\{(q_i,\omega_i)\}_{i=1}^N$ is the weighted point cloud with
$q_i=(y_i,-b_i)$ and 
$\omega_i=-2\,b_i^2$, and where 
$\Sigma_i^+$ is one sheet of a hyperboloid given by $\Sigma_i^+=\left\{\(x, |X-Y_i|+b_i\):X=(x,\alpha), x\in \R^2\right\}$.
Therefore,
$$
\Lag_i(\b)=\proj_{\R^2\times\{\alpha\}}\left( 
\Pow_i(Q) \cap \Sigma_i^+\cap \Omega
\right).
$$
\end{proposition}

In practice, condition \eqref{eq:bi and bj are not too far} is not restrictive because to use the damped Newton algorithm from Section \ref{sec:algorithm}, one needs to assume that the Laguerre cells at the initialization vector $\b^0$ are non empty which by the corollary below implies \eqref{eq:bi and bj are not too far}.
\begin{corollary}
If the vector $\b$ satisfies $\Lag_i(\b)\neq \emptyset$ for each $1\leq i \leq N$, then
$$
\Lag_i(\b)=\proj_{\R^2\times\{\alpha\}}\left( 
\Pow_i(Q) \cap \Sigma_i^+\cap \Omega
\right),
$$
for all $1\leq i \leq N$.
\end{corollary}
\begin{proof}  Fix $1\leq i,j\leq N$ with $i\neq j$ and suppose that $\b=(b_1,\cdots ,b_N)$ is such that $\Lag_i(\b)\neq \emptyset$ and $\Lag_j(\b)\neq \emptyset$ contain points $X_i$ and $X_j$ respectively. Then,
  $$\nr{X_i} + \nr{X_i - Y_i} =  c(X_i,Y_i) + b_i \leq c(X_i,Y_j) + b_j = \nr{X_i} + \nr{X_i - Y_j} + b_j. $$
  Thus, $b_i - b_j \leq |Y_i-Y_j|=|y_i-y_j|$. By symmetry, we get
  \begin{equation*}\label{eq:bi and bj are not too far}
 |b_i - b_j| \leq \nr{y_i - y_j} < \sqrt{4(\alpha-\beta)^2+|y_i-y_j|^2}.\qedhere
\end{equation*}
  \end{proof}

\begin{remark}\label{rm:initialization vector}\rm
If $\alpha=\beta$, then the Laguerre diagram is known as \emph{Apollonius diagram} where all the surfaces $\Sigma_i^+$ are half-cones. When $\alpha\neq\beta$, then each $\Sigma_i^+$ is not anymore a cone but a sheet of a hyperboloid. 
\end{remark}
\begin{remark}[Initialization]\rm
To apply the algorithm from Section \ref{sec:algorithm}, we need to find an initial vector $\b^0=(b_1,\cdots,b_N)$ for which the corresponding Laguerre cells are not empty. Taking $\b^0=0$ this is the case when 
$$
  \proj_{\R^2\times\{\alpha\}}(Y_i) \in \Omega, \forall i \in \{1,\cdots, N\}.
$$
Indeed, if we denote by $X_i=(y_i,\alpha)$ such a point, we see that $|X_i-Y_i| = |\beta-\alpha| \leq |X_i-Y_j|$ for any $j$, implying that $X_i \in \Lag_i(0)$.
For other options to choose the initialization vector see \cite[Sec. 2.3]{2019-meyron-initialization-procedures}.
\end{remark}

\subsection{Proof of Proposition~\ref{prop:computationlaguerre_PSNF}}
Let us consider the hyperboloid of two sheets
\[
\Sigma_i=\left\{(x,x_3): |X-Y_i|^2=(x_3-b_i)^2, X=(x,\alpha), x\in \R^2\right\}
\]
with $Y_i=(y_i,\beta)$, $\beta>\alpha$.
The upper sheet of this hyperboloid is given parametrically by
\[
\Sigma_i^+=\left\{\(x,  |X-Y_i|+b_i\):X=(x,\alpha), x\in \R^2\right\}
\]
and the lower sheet is given by
\[
\Sigma_i^-=\left\{\(x, - |X-Y_i|+b_i\):X=(x,\alpha), x\in \R^2\right\}.
\]
Clearly, $\Sigma_i^-$ and  $\Sigma_i^+$ are symmetric with respect to the hyperplane $\{x_3=b_i\}$, 
and $\Sigma_i = \Sigma_i^- \cup \Sigma_i^+$.

We first need to determine the relative positions of the hyperboloids $\Sigma_i$ and $\Sigma_j$.
\begin{lemma}\label{lemma:above}
We have the following:
\begin{enumerate}
\item If $\Sigma_i^+ \cap \Sigma_j^- = \emptyset$, then $\Sigma_i^+$ is strictly above $\Sigma_j^-$, i.e., for any $x\in \R^2$ such that $(x,x_3^i)\in \Sigma_i^+$ and $(x,x_3^j)\in \Sigma_j^-$, one has $x_3^i > x_3^j$.
\item  
$
\sqrt{4(\alpha-\beta)^2+|y_i-y_j|^2} > b_j - b_i$
 if and only if 
$\Sigma_i^+ \cap \Sigma_j^- = \emptyset$.
\item If $\Sigma_i^+ \cap \Sigma_j^- \neq \emptyset$, then $\Sigma_i^+$ is strictly below $\Sigma_j^+$, i.e., for any $x\in \R^2$ such that $(x,x_3^i)\in \Sigma_i^+$ and $(x,x_3^j)\in \Sigma_j^+$, one has $x_3^i < x_3^j$.   
\end{enumerate}
\end{lemma}
\begin{proof}
(1) Since $\Sigma_i^+$ opens upwards and $\Sigma_j^-$ opens downwards, if they don't intersect, it is clear that  $\Sigma_i^+$ must be strictly above $\Sigma_j-$.

(2) We have $\Sigma_i^+ \cap \Sigma_j^- \neq \emptyset$ iff there exists $X\in \R^2\times\{\alpha\}$ with
$$
|X-Y_i| + |Y_j-X| = b_j -b_i.
$$
On the other hand, for each $X\in \R^2\times\{\alpha\}$, one has
$$
|X-Y_i| + |Y_j-X| \geq |X_{i,j}-Y_i| + |X_{i,j} - Y_j| = 2\sqrt{(\alpha-\beta)^2+|y_i-y_j|^2/4},
$$
where $X_{i,j}=((y_i+y_j)/2,\alpha)$. Therefore, $\sqrt{4(\alpha-\beta)^2+|y_i-y_j|^2} > b_j - b_i$ implies  $\Sigma_i^+ \cap \Sigma_j^- = \emptyset$. 

Vice versa, from Item (1), $|X-Y_i|+b_i>-|X-y_j|+b_j$ for each $X=(x,\alpha)$ and so $b_j-b_i$ satisfies the desired inequality.

(3) By contradiction. Suppose there are points $(x,x_3^i) \in \Sigma_i^+$ and $(x,x_3^j) \in \Sigma_j^+$ 
with $x_3^i\geq x_3^j$.
If $X =(x,\alpha)$, then
$
|X-Y_i| + b_i =x_3^i \geq x^3_j = |X-Y_j| + b_j,
$
and so $b_j-b_i\leq |Y_i-Y_j|$ from triangle inequality.
Now, $\sqrt{4(\alpha-\beta)^2+|y_i-y_j|^2} \leq  b_j - b_i$ from (2). 
Since $|Y_i-Y_j|<\sqrt{4(\alpha-\beta)^2+|y_i-y_j|^2}$, we obtain a contradiction.

\end{proof}

The previous lemma leads to the following.

\begin{lemma}\label{lm:projections of Lij}
Fix $i,j$ and define the sets 
$$
L_{i,j}^\leq := \left\{ X\in \R^2\times\{\alpha\}: |X-Y_i| + b_i \leq |X-Y_j|+b_j\right\} \supset \Lag_i(\b),
$$
$$
H_{ij}^\leq := \left\{x\in \R^3: |x - q_i|^2 + \omega_i \leq |x - q_j|^2 + \omega_j  \right\},\quad \omega_j=-2\,b_j^2; \quad q_j=(y_j,-b_j).
$$
We have:
\begin{enumerate}
\item If $\sqrt{4(\alpha-\beta)^2+|y_i-y_j|^2} > b_j - b_i$, then
$$
L_{i,j}^{\leq}= \proj_{\R^2\times\{\alpha\}}\left( H_{ij}^\leq \cap \Sigma_i^+ \right)
$$
\item If $\sqrt{4(\alpha-\beta)^2+|y_i-y_j|^2} \leq b_j - b_i$, then $L_{i,j}^{\leq}=\R^2\times\{\alpha\}$.

\end{enumerate}
\end{lemma}
\begin{proof} 
(1) Let $X=(x,\alpha)\in L_{i,j}^\leq$, and $x_3=|X-Y_i|+b_i$. 
By Lemma~\ref{lemma:above} (2) and (1), the point $(x,x_3)\in \Sigma_i^+$ is above $\Sigma_j^-$, and so $x_3 - b_j > -|Y_j-X|$. 
By definition of $L_{i,j}^\leq$, we also have $x_3 - b_j \leq |Y_j-X|$, which implies 
$$
(x_3 - b_i)^2 = |X-Y_i |^2 \quad  \mbox{and} \quad   (x_3 - b_j)^2 \leq  |X-Y_j|^2.
$$
Expanding these two equations, one gets 
$$
\left\{
\begin{array}{l}x_3^2 - 2 b_ix_3 +b_i^2 = |X|^2-2 X\cdot Y_i+ |Y_i|^2\\
x_3^2 - 2 b_jx_3 +b_j^2 \leq |X|^2-2 X\cdot Y_j+ |Y_j|^2.\\
\end{array}\right.\\
$$
Subtracting the first line from the second line yields
$$
-2(x,x_3) \cdot (y_i, -b_i) + |Y_i|^2 - b_i^2 
\leq -2(x,x_3) \cdot (y_j, -b_j) + |Y_j|^2 - b_j^2
$$
which can be rewritten as 
\[
|(x,x_3) - (y_i,-b_i)|^2 - 2 b_i^2 \leq  |(x,x_3) - (y_j,-b_j)|^2 - 2 b_j^2. 
\]
This means $(x,x_3) \in H_{ij}^\leq$ and so $(x,x_3) \in H_{ij}^\leq \cap \Sigma_i^+$.
To show the opposite inclusion, let $(x,x_3)\in  H_{ij}^\leq \cap \Sigma_i^+$ and put $X=(x,\alpha)$. Then one has 
\[
\left\{
\begin{array}{l}
x_3=|X-Y_i|+b_i\\
|(x,x_3) - (y_i,-b_i)|^2 - 2b_i^2 \leq  |(x,x_3) - (y_j,-b_j)|^2 - 2 b_j^2.
\end{array}\right.\\
\]
Reversing the previous calculation, one gets $(x_3 - b_j)^2 \leq  |X-Y_j|^2$. This obviously implies $-|X-Y_j|\leq x_3 - b_j \leq  |X-Y_j|$, which gives in particular that  $X\in L_{i,j}^\leq$, completing the proof of (1).

(2) From Lemma~\ref{lemma:above} (2), $\Sigma_i^+ \cap \Sigma_j^- \neq \emptyset$ and so from Item (3) in that lemma, $\Sigma_i^+$ is strictly below $\Sigma_j^+$. That is, $b_i+|X-Y_i|<b_j+|X-Y_j|$ for all $x\in \R^2$, $X=(x,\alpha)$. This means 
$L_{i,j}^\leq$ is the whole plane $\R^2\times\{\alpha\}$.
\end{proof}

\begin{proof}[Proof of Proposition~\ref{prop:computationlaguerre_PSNF}]
Let $i\in \{1,\cdots, N\}$. From \eqref{eq:bi and bj are not too far} we can apply 
Lemma \ref{lm:projections of Lij}(1) to obtain 
\[
H_i(\b) = 
\displaystyle\bigcap_{1\leq j \leq N} L_{i,j}^{\leq} = 
\bigcap_{1\leq j \leq N} \proj_{\R^2\times\{\alpha\}}\left( H_{ij}^\leq \cap \Sigma_i^+ \right)
= \proj_{\R^2\times\{\alpha\}} \left( \Sigma_i^+ \cap \bigcap_{1\leq j \leq N}   H_{ij}^\leq \right)
\]
Since by definition $\Pow_i(Q) = \bigcap_{1\leq j \leq N}   H_{ij}^\leq$,
the proposition follows.
\end{proof}

\subsection{Numerical experiments} 
In all three numerical experiments, we assume that the source measure
$\rho\equiv \frac{1}{4}$ is uniform over the square $\Omega =
[-1,1]^2\times \{\alpha\}$, and we also assume that the metasurface is at height $\alpha = 1$. 
The Laguerre cells are computed using
Proposition~\ref{prop:computationlaguerre_PSNF}, by intersecting 3D
power cells with a quadric. To describe the algorithm, we use the
notation from that proposition:
\begin{itemize}
\item
First, we compute 3D the power diagram $(\Pow_i(Q))_{1\leq i\leq N}$
using the CGAL library and we restrict each cell to the lifted domain
$\Omega\times \R$ by computing the intersection $P_i = \Pow_i(Q)\cap
(\Omega\times \R)$.  In the implementation, we assume that
$\Omega\subseteq\R^2$ is a convex polygon.
\item For every pair $i\neq j \in \{1,\hdots,N\}$, we compute the curve
$\gamma_{ij}$ corresponding to the projection on $\R^2\times \{\alpha\}$
of the facet $P_i\cap P_j$ with the quadric
$\Sigma_i^+$,
$$\gamma_{ij} = \proj_{\R^2\times\{\alpha\}}(P_i \cap
P_j \cap \Sigma_i^+).$$ In practice, we need to make this computation
only if the power cells already have a non-empty interface, i.e. if
$P_i\cap P_j \neq \emptyset$. We also note that by
Proposition~\ref{prop:computationlaguerre_PSNF}, 
$$ \Lag_i(\b) \cap \Lag_j(\b) = \gamma_{ij}. $$
\item We finally compute the intersection of each Laguerre cell with the
boundary of the domain using the formula
$$\gamma_{i,\infty} = \Lag_i(\b)\cap\Omega
= \proj_{\R^2\times\{\alpha\}}(\Pow_i(Q)\cap \partial \Omega).$$
\end{itemize}
By construction, the boundary of the $i$th Laguerre cell is given by
$$ \partial \Lag_i(\b) = \gamma_{i,\infty} \cup \bigcup_{j\neq
i} \gamma_{ij}, $$ and the union is disjoint up to a finite set, which
is negligible. We may then use this description of the boundary of the
cell $\Lag_i(\b)$ to compute the integral of $\rho$ over $\Lag_i(\b)$
using divergence theorem 
$$ G_i(\b) = \int_{\Lag_i(\b)} \rho(X)d X
= \frac{1}{8} \int_{\partial \Lag_i(\b)} X \cdot n_{i}(X)
dX, $$
where $n_i$ denotes the exterior normal to $\Lag_i(\b)$, and where the second integral is 1-D.
The partial derivatives are computed
using \eqref{eq:derivatives of Gi with respect bj j not equal i}: for
$i\neq j$ we have
$$\dfrac{\partial G_i}{\partial b_j}(\b)
= 
\int_{\gamma_{ij}} \dfrac{1}{4}\dfrac{1}{\left| \nabla_X c(X,Y_i)-\nabla_X c(X,Y_j)\right|}\,dX,$$
where again the integrals are one dimensional.
The code to compute the intersection between the power cell
$\Pow_i(Q)$ and the quadric $\Sigma_i^+$ and to perform the numerical
integration is written in a combination of C++ and Python, and is
available online\footnote{\url{https://github.com/mrgt/ot-optics}}, as
well as the experiments presented below.

\subsubsection{Effect of changes in $\beta-\alpha$ on the shape of $\Lag_i(\b)$}
\label{subsec:numerics-beta}
In the first numerical experiment, we study the effect of the distance
between the metasurface and the target, $\delta = \beta-\alpha$, on
the shape of the Laguerre cells. We assume that the target is of the
form
$$ \nu = \frac{1}{N} \sum_{1\leq i \leq 25} \delta_{y_i}, $$ where
$N=25$ and $\{y_1,\hdots,y_{N}\}$ is a uniform $5\times 5$ grid
contained in the square $[0,1]^2$, and we solve the optimal transport
problem between $\rho$ and $\nu$. Our goal in this first experiment is
to visualize the effect of changes in $\delta$, the vertical
distance between the source and the metasurface, on the shape of the
solution. We initialize the damped Newton algorithm described
in \Cref{sec:algorithm} with $\b^0 = (0,\hdots,0)$. The associated
Laguerre cells coincides with the Voronoi cells of the point cloud
$\{Y_1,\hdots,Y_N\}$, i.e.
$$\mathrm{Vor}_i = \{ X \in  \R^2\times \{\alpha\} \mid \forall
j\in\{1,\hdots,N\},~ |{X - Y_i}| \leq |{X - Y_j}|\}.$$ and is shown
on the first row and column of \Cref{fig:numerics-beta}. We solve the
near-field metasurface problem for several values of $\delta$, and we
display the Laguerre cells of the solution. In particular, one can see
that for $\delta = 2$, the solution is very similar to the solution of
the optimal transport problem for the ``standard'' quadratic cost.

\subsubsection{Convergence speed} \label{subsec:numerics-gaussian}
\begin{figure}[H]
\includegraphics[width=0.45\textwidth]{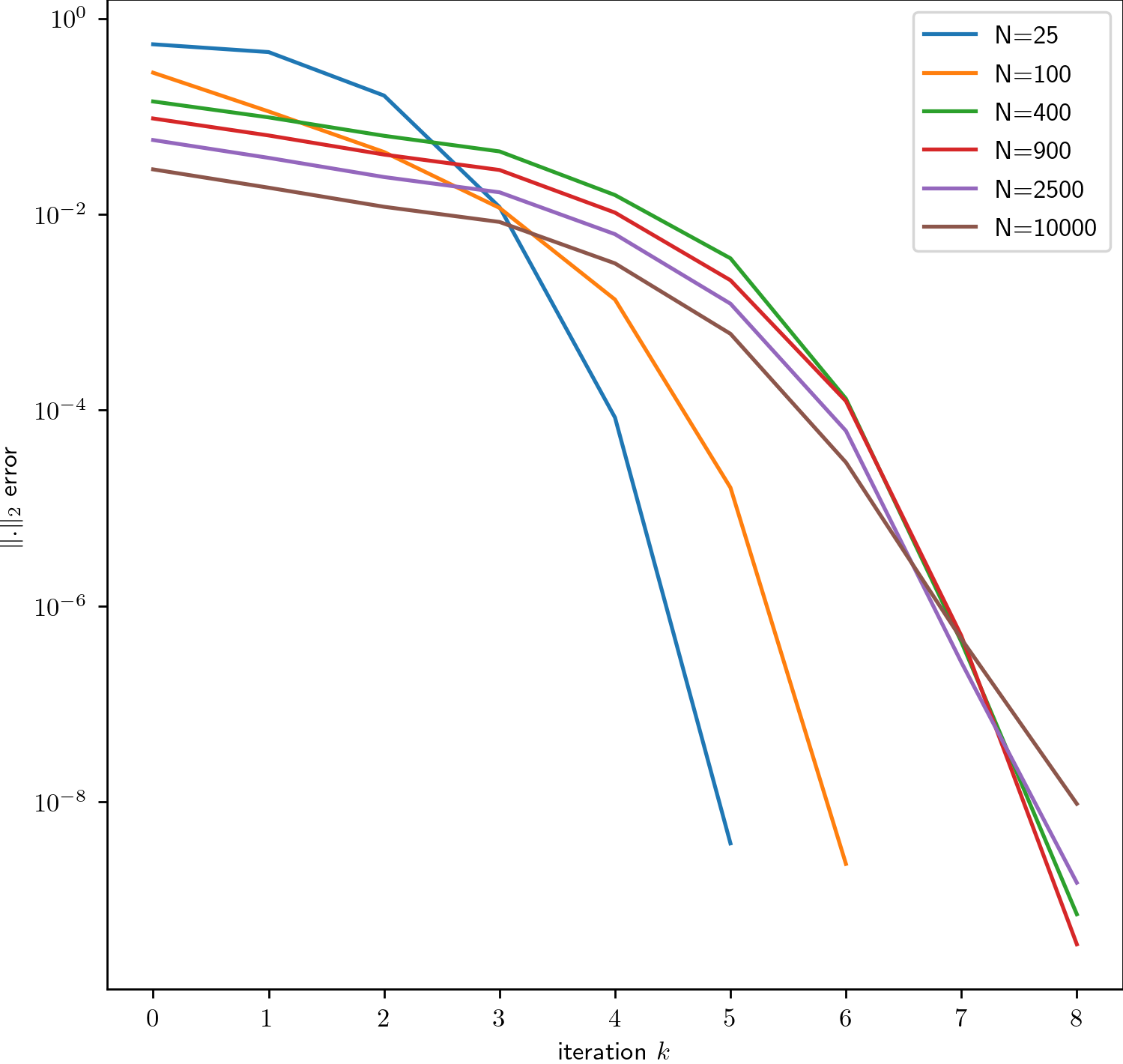}
    \caption{Convergence of the numerical error in terms of the iteration number for the
    near-field metasurface problem described 
    in \S\ref{subsec:numerics-gaussian}, for
    $N\in\{5,10,20,30,40,50,100\}^2$.  \label{fig:numerics-gaussian}}
\end{figure}

\begin{figure}[H]
\includegraphics[width=0.5\textwidth]{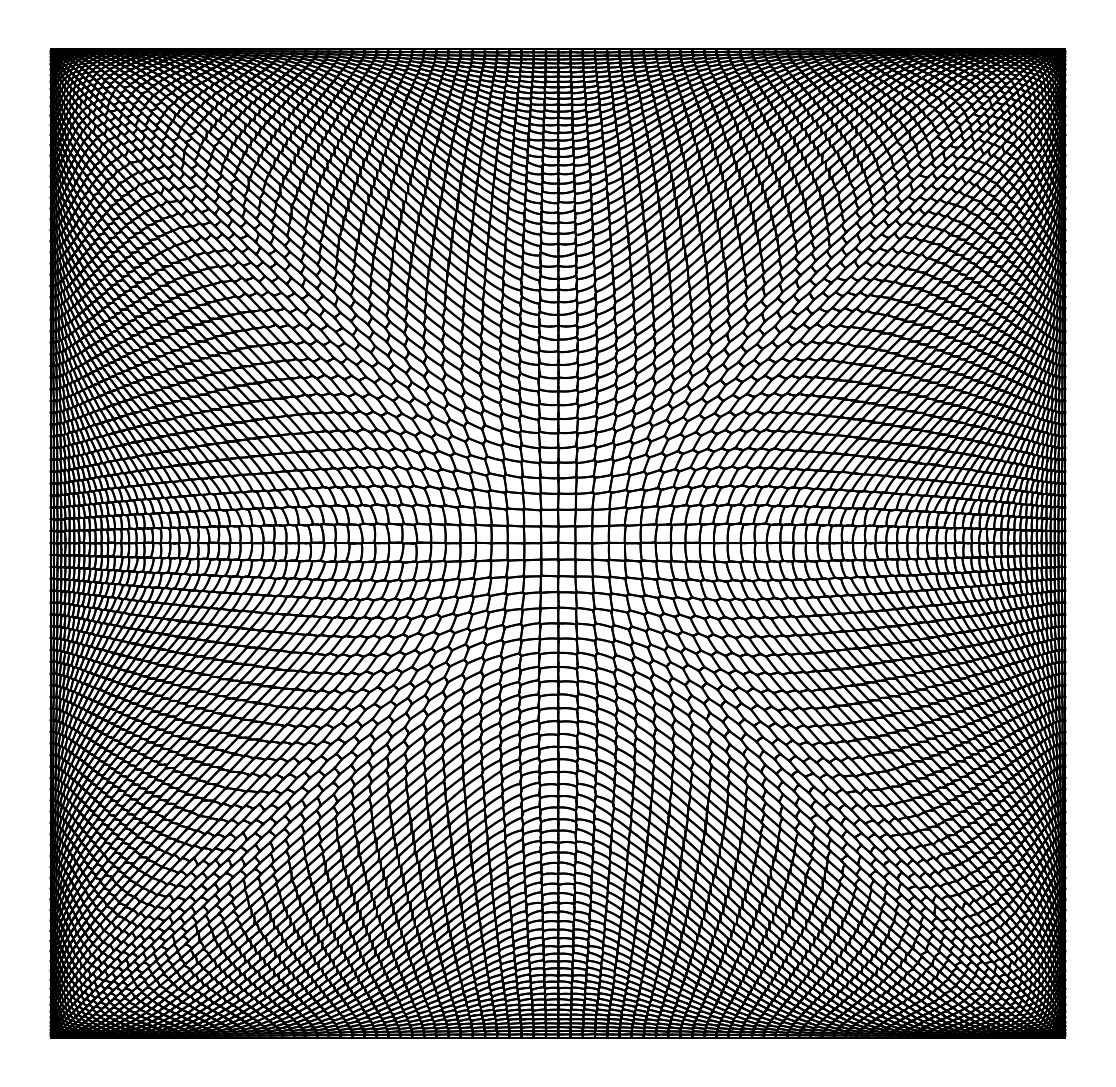}
    \caption{Laguerre diagram associated to the solution of the
     near-field metasurface problem  described
    in \S\ref{subsec:numerics-gaussian}, for
    $N=100^2$.  \label{fig:numerics-gaussian2}}
\end{figure}

In our second numerical experiment, the target measure $\nu$
approximates the restriction of the Gaussian $e^{-2|{\cdot}|^2}$ to
the unit square $[-1,1]^2\times \{\alpha\}$. More precisely, the measure $\nu$ is of the
form
$$ \nu = \frac{1}{N} \sum_{1\leq i\leq N} \nu_i \delta_{y_i}, $$
where $N = n^2$ and $n\in \{5,10,20,30,40,50,100\}$. The points $\{y_1,\hdots,y_N\}$ form a uniform $n\times n$ grid in the square $[-1,1]^2\times\{\alpha\}$. The mass $\nu_i$ of the Dirac $\delta_{y_i}$ is defined by evaluating  a Gaussian at $y_i$:
$$ \nu_i = e^{-2|{y_i}|^2}/\sum_{1\leq j\leq N}
e^{-2|{y_j}|^2}. $$ Figure \ref{fig:numerics-gaussian2} displays the
solution of this problem for
$N=100^2$. Figure \ref{fig:numerics-gaussian} displays the decrease
of the  numerical error along the iterations of the algorithm, defined as
$$ \varepsilon_k = \left(\sum_i (H_i(\b^k) - \nu_i)^2\right)^{1/2}, $$
for several values of $N$. In particular, one can see from this figure that a numerical error of $10^{-8}$ is reached in less than 8 iterations, even for $N=10^4$.

\subsubsection{Visualization of the phase} \label{subsec:phasevisu}
\begin{figure}[H]
  \begin{center}
    \begin{tabular}{ccc}
\includegraphics[width=0.27\textwidth]{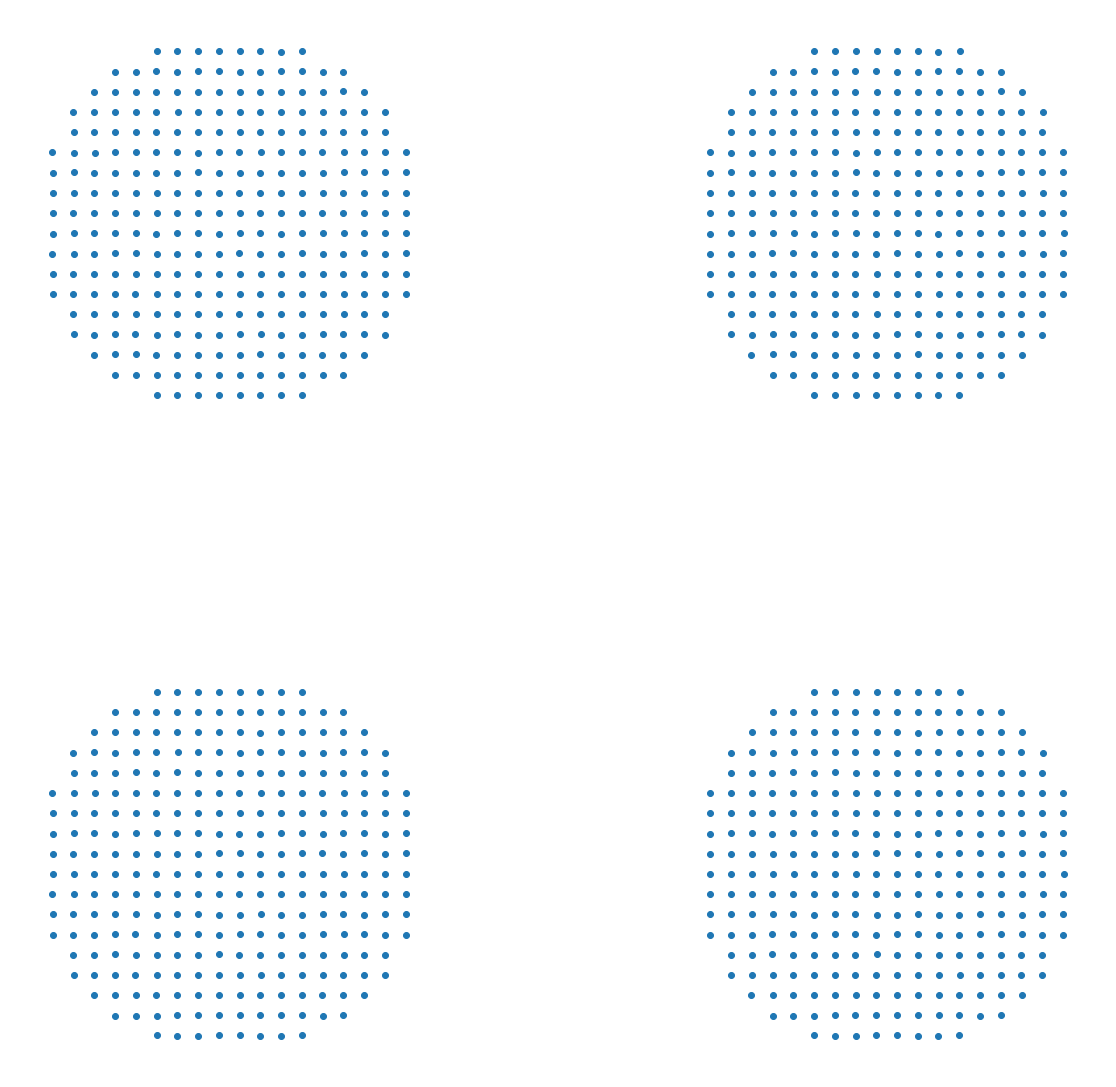}&
\includegraphics[width=0.27\textwidth]{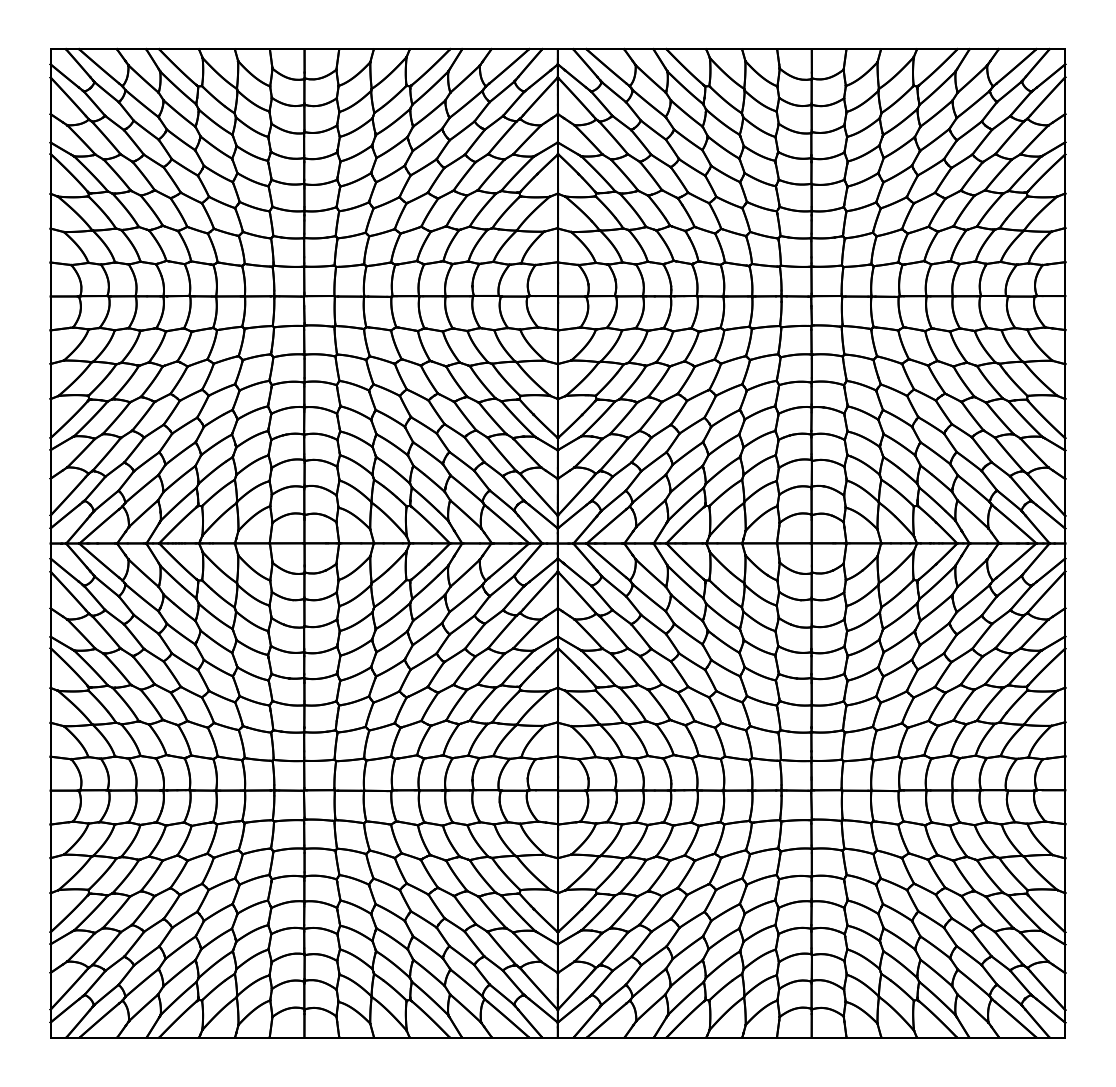}&
\includegraphics[width=0.27\textwidth]{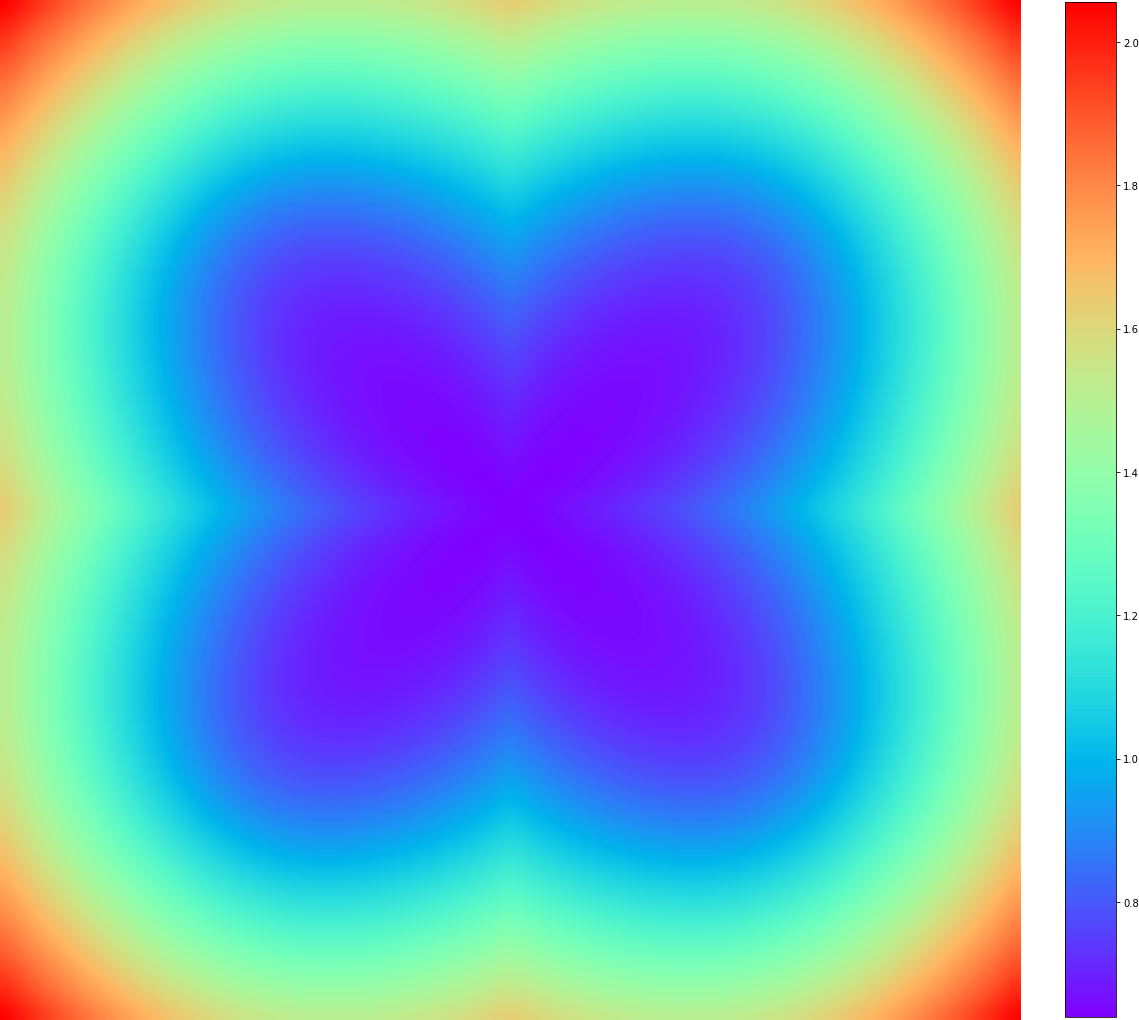}\\
\includegraphics[width=0.27\textwidth]{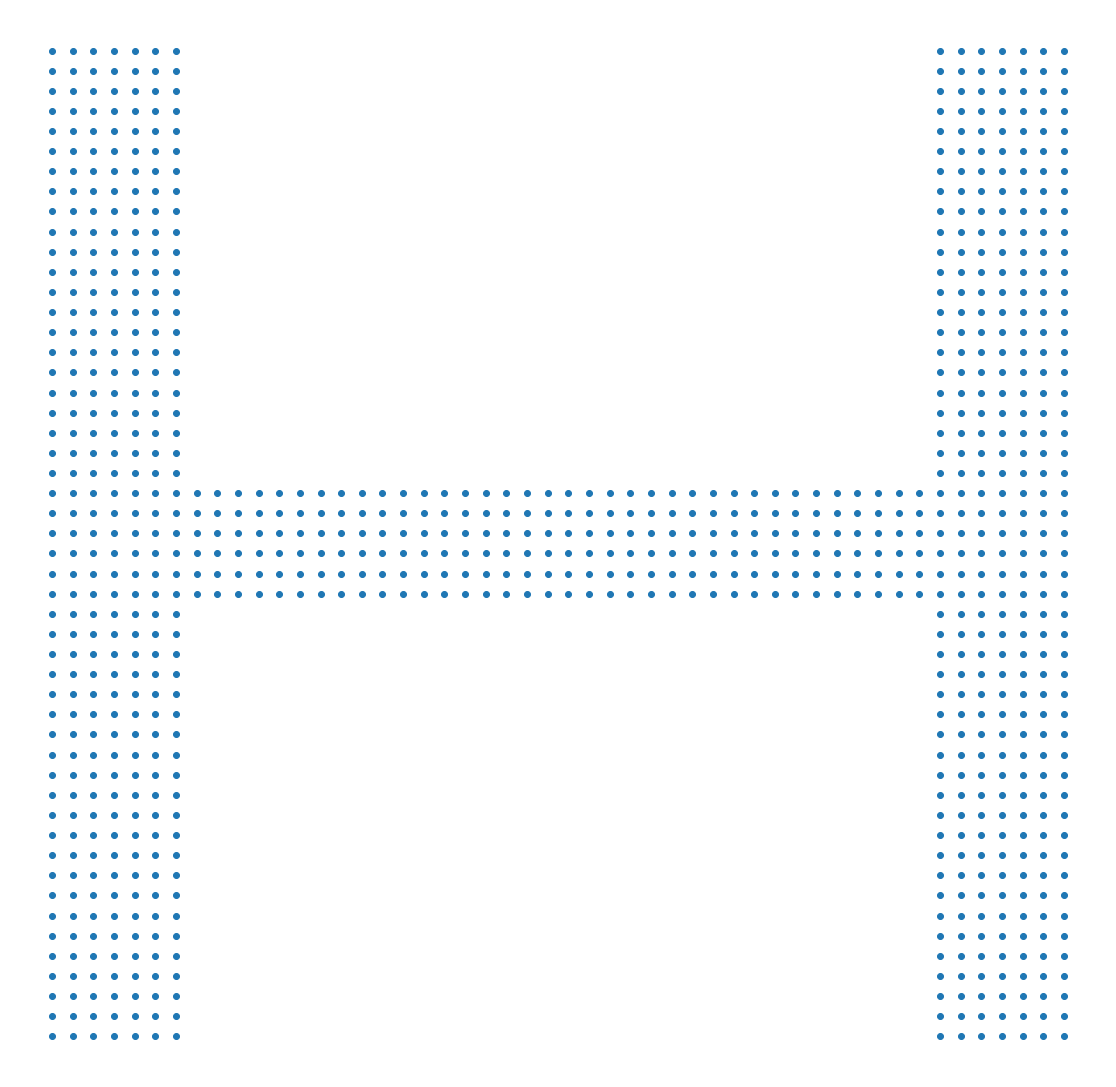}&
\includegraphics[width=0.27\textwidth]{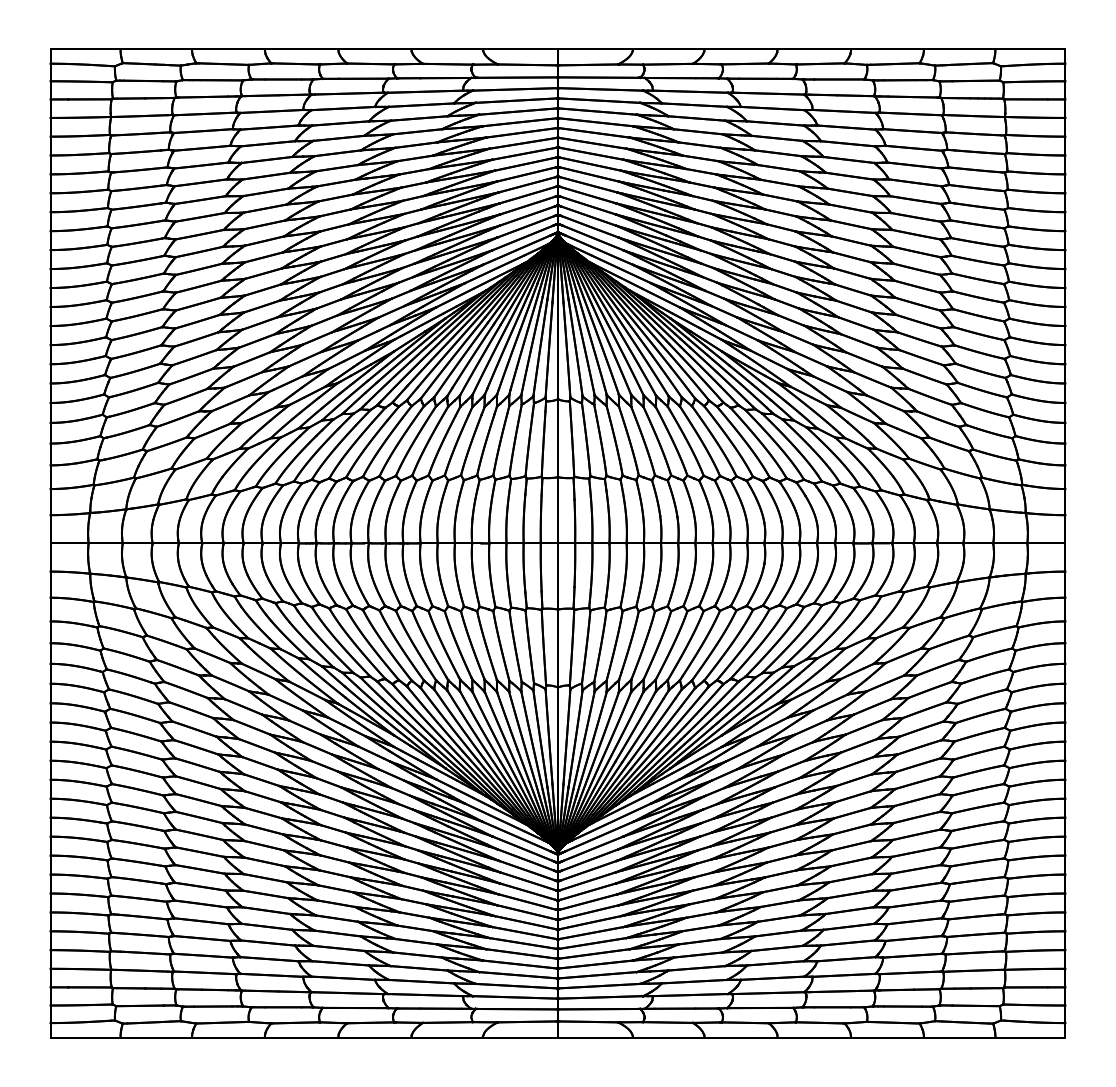}&
\includegraphics[width=0.27\textwidth]{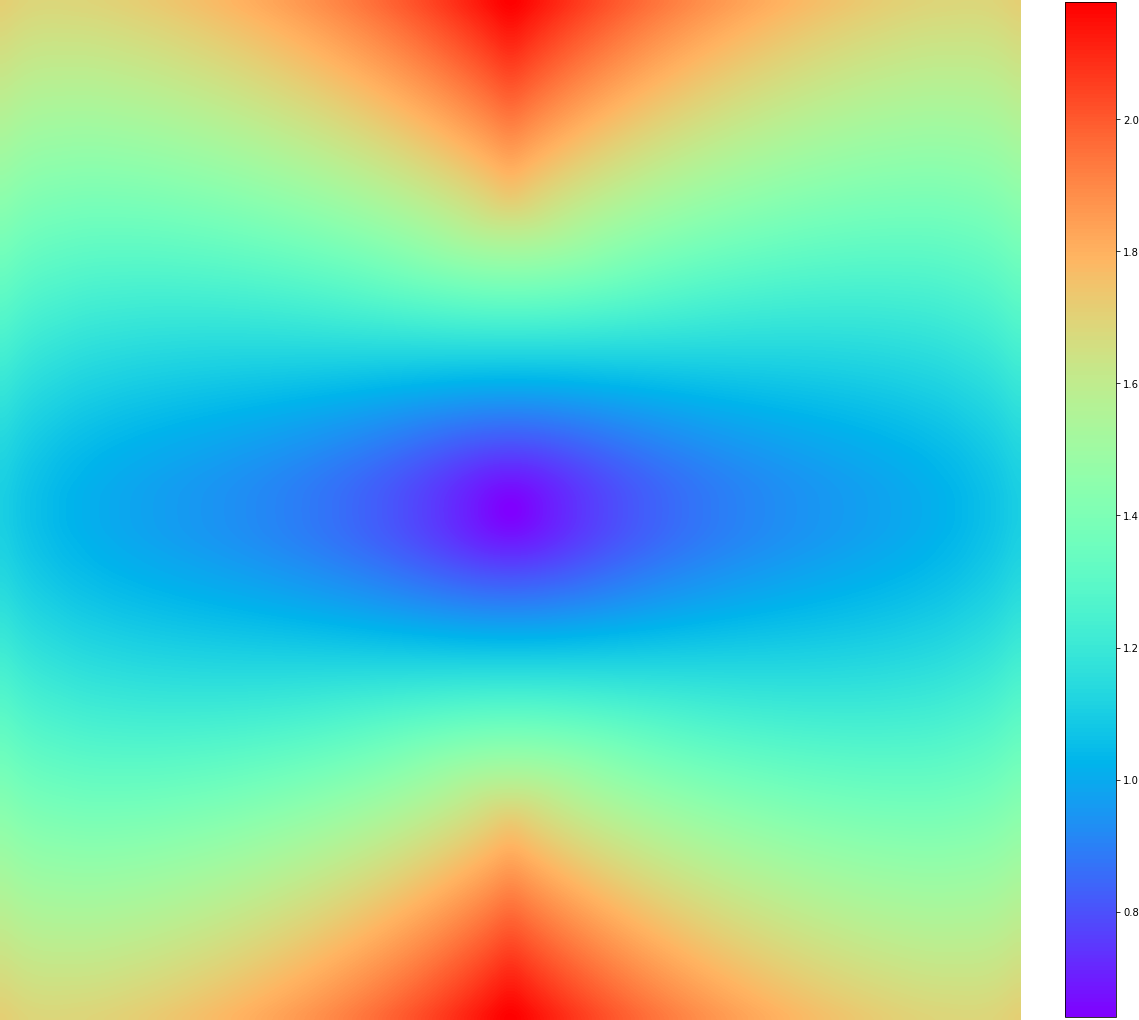}
\end{tabular}
\caption{Near-field metasurface with a target measure supported over four discretized disks (top) or a discretized letter H (bottom): Support of the target measure (left), Laguerre cells associated to the solution (middle), and the  corresponding phase discontinuity $\phi$ (right).
\label{fig:4disks}}
\end{center}
\end{figure}

In this last numerical experiment, the target is uniform over four
discretized disks (Figure~\ref{fig:4disks}, top row) or over a
discretized letter H (Figure~\ref{fig:4disks}, bottom row), i.e. $\nu
\equiv 4/N$ where $N$ is the number of points composing the
discretized shapes. Figure~\ref{fig:4disks} displays the Laguerre
cells corresponding to the solution of the near-field metasurface
problem. We also display the corresponding phase discontinuity $\phi$,
which can be computed thanks to Theorem~\ref{thm:existence of
  metasurfaces near field}. On the ``four disks'' example, one may
notice that the gradient of the phase discontinuity $\phi$ seems to
exhibit a discontinuity on the ``cross'' $\{0\}\times
[0,1]\times\{\alpha\}\cup [0,1]\times\{0\}\times\{\alpha\}$: this
corresponds to a jump in the transport map which is necessary to cross
the void between the four disks. 

%% file: section5-farfield.tex
\section{Refraction into the far field}\label{sec:refraction into the far field}

We consider here the case when incident rays emanate in a collimated beam and the case when they emanate from a point source. Since the arguments to treat these problems are similar to the ones used in Section \ref{sec:existence and uniqueness}, we only indicate the modifications that are needed.
\subsection{Collimated beam}
Let $\Gamma$ be the horizontal plane $x_3=a$ in $\R^3$, $\e=(0,0,1)$.
The phase discontinuity function $\phi(x)$ with $x=(x_1,x_2,x_3)$ defined in a neighborhood of $\Gamma$ such that the metasurface $(\Gamma,\phi)$ refracts all vertical rays having direction $\e$ into a fixed unit direction ${\bf m}=\(m_1,m_2,m_3\)$, with $m_3>0$, is given by 
\[
\phi(x)={\bf v}\cdot x,\text{ with } {\bf v}=(\m\cdot \e)\,\e-\m=\(-m_1,-m_2,0\);
\]  
see \cite[Theorem 4.1]{gutierrezsabra:chromaticaberrationinmetasurfaces} for a proof in the more general case when $\m$ is a variable set of directions depending on $x$.
Notice that if ${\bf m}=\(m_1,m_2,m_3\)$ and ${\bf m}'=\(m_1',m_2',m_3'\)$ are two unit vectors with $m_3,m_3'>0$ and $\(m_1,m_2\)=\(m_1',m_2'\)$, then ${\bf m}={\bf m}'$.

Each vector ${\bf m}=\(m_1,m_2,m_3\)\in S^2$ with $m_3>0$ can be identified with its projection $\(m_1,m_2\)$ on the disk of radius one around $(0,0)$ where $m_3=\sqrt{1-m_1^2-m_2^2}$.

Fix a compact region $\Omega\subset \Gamma$ with points $x=(x_1,x_2,a)$ and a compact region $\Omega^*\subset S^2_+$.

\begin{definition}[Admissible phase]
The function $\phi:\Omega\to \R$ is an admissible phase refracting $\Omega$ into $\Omega^*$ if for each $x_0\in \Omega$ there exists ${\bf m}\in \Omega^*$ and $b\in \R$  such that 
\begin{align*}
\phi(x)&\geq b+\((\m\cdot \e)\,\e-\m\)\cdot x:=L(\m,b,x),\qquad \forall x\in \Omega,\text{ and}\\
\phi(x_0)&=b+\((\m\cdot \e)\,\e-\m\)\cdot x_0.
\end{align*}
We say that $L(\m,b,x)$ is a supporting phase to $\phi$ at $x_0$.
\end{definition}

\begin{definition}
Given $x_0\in \Omega$ and $\phi$ an admissible phase,
we define the set-valued mapping $\mathcal N_\phi:\Omega\to \mathcal P(\Omega^*)$ given by
\[
\mathcal N_\phi (x_0)=\{m\in \Omega^*:\text{there exists $b\in \R$ such that $L(\m,b,x)$ is a supporting phase to $\phi$ at $x_0$} \};
\]
and the inverse map
\[
\(\mathcal N_\phi\)^{-1}(m)=\{x\in \Omega:m\in \mathcal N_\phi (x)\}.
\] 
\end{definition}
Define the class
\[
\mathcal F=\{\phi:\Omega\to \R: \text{$\phi$ is an admissible phase refracting $\Omega$ into $\Omega^*$}\}.
\]
Given $f\in L^1(\Omega)$, non negative, assuming that $\partial \Omega$ has 2-dimensional Lebesgue measure zero, and given 
a Borel measure $\mu$ in $\Omega^*$ satisfying 
$
\int_{\Omega}f(x)\,dx=\mu(\Omega^*),
$
our problem is then to find $\phi\in \mathcal F$ solving
\begin{equation}\label{eq:definition of phase solution}
\mathcal M_\phi(E):=\int_{(\mathcal N_\phi)^{-1}(E)}f(x)\,dx=\mu(E)
\end{equation}
for each Borel set $E\subset \Omega^*$.
To show existence and uniqueness to this problem, we proceed as in Section \ref{sec:existence and uniqueness} with the following changes.
The class $\mathcal F$ satisfies the following properties that follow immediately from the definitions above
\begin{enumerate}
\item[(A1)'] If $\phi_1,\phi_2\in \mathcal F$, then $\phi_1\vee\phi_2=\max\{\phi_1,\phi_2\}\in \mathcal F$,
\item[(A2)'] if $\phi_1(x_0)\geq \phi_2(x_0)$, then $\mathcal N_{\phi_1}(x_0)\subset \mathcal N_{\phi_1\vee \phi_2}(x_0)$,
\item[(A3)'] Given $\m\in \Omega^*$ the functions $L(\m,b,x)\in \mathcal F$, $b\in (-\infty,\infty)$,  satisfy the following
\begin{enumerate}
\item[(a)] $\m\in \mathcal N_{L(\m,b,x)}(x)$ for all $x\in \Omega$,
\item[(b)] $L(\m,b,x)\leq L(\m,b',x)$ for all $b\leq b'$,
\item[(c)] for each $\m\in \Omega^*$, $L(\m,b,x)\to +\infty$ uniformly for $x\in \Omega$ as $b\to +\infty$,
\item[(d)] for each $\m\in \Omega^*$, $\max_{x\in \Omega} |L(\m,b',x)-L(\m,b,x)|\to 0$ as $b'\to b$.
\end{enumerate}
\end{enumerate}
Recalling the notation at the beginning of Section \ref{sec:existence of solutions},
we now let $\X=\Omega$ and $\Y=\Omega^*$, and with similar arguments but now using conditions $(A1)'-(A3)'$ instead of $(A1), (A2), (A3)''$,  
%
%
we get that
$
\mathcal N_\phi\in C_s(\Omega,\Omega^*)
$
for each $\phi\in \mathcal F$.
Then to prove existence and uniqueness when $\mu=\sum_{i=1}^{N} g_{i}\,\delta_{{\bf m}_i}$, we use \cite[Theorem 2.12]{gutierrez-huang:nearfieldrefractor},
for which we only need to verify its hypotheses.
In fact, we need to verify that there exist numbers $b_1^0,\cdots ,b_N^0$ such that 
the admissible phase $\phi_0(x)=\max_{1\leq i\leq N}L\({\bf m}_i,b_i^0,x\)$, for $x\in \Omega$, satisfies
$\mathcal M_{\phi_0}({\bf m}_j)\leq g_j$ for $2\leq j\leq N$.
By continuity, given any $b_1^0$ we can pick $b_2^0,\cdots ,b_N^0$ tending to $-\infty$ such that 
$L\({\bf m}_i,b_i^0,x\)< L\({\bf m}_1,b_1^0,x\)$ for all $i\neq 1$ and $x\in\Omega$.
Therefore, $\phi_0(x)=L\({\bf m}_1,b_1^0,x\)$.
Since the points ${\bf m}_i=\(m_1^i,m_2^i,m_3^i\)\neq {\bf m}_j=\(m_1^j,m_2^j,m_3^j\)$ for $i\neq j$, the family of planes having equations $z=\alpha-m_1^i\,x_1-m_2^i\,x_2 $ are never parallel. 
Consequently, $\mathcal N_{\phi_0}({\bf m}_i)\subset \partial \Omega$ for $i\neq 1$ and so
$\mathcal M_{\phi_0}({\bf m}_i)=0$ for $i\neq 1$. 
The hypotheses in \cite[Theorem 2.12]{gutierrez-huang:nearfieldrefractor} then hold in our case.
In addition, from the conditions $(A1)'-(A3)'$ above we can apply \cite[Theorem 2.12]{gutierrez-huang:nearfieldrefractor} to obtain the following.
\begin{theorem}\label{thm:existencesumofdiracmassessecondbdryMA}
Let $\m_{1},\cdots , \m_{N}$ be distinct points in $\Omega^*$,  
$g_{1},\cdots , g_{N}$ are positive numbers, and $f\in L^{1}(\Omega)$
with
\begin{equation}\label{eq:conservationofenergysecondbdryMA}
\int_{\overline{\Omega}} f(x)\,dx=\sum_{i=1}^{N} g_{i};
\end{equation}
$\mu=\sum_1^N g_i\,\delta_{\m_i}$.
Then 
given any $b_1\in \R$, there exist numbers $b_{2},\cdots , b_{N}$ such that the convex function
\begin{equation}\label{eq:solutionsecondbdryMA}
\phi(x)=\max_{1\leq i \leq N}\{ L(\m_i,b_i,x)\}
\end{equation}
solves \eqref{eq:definition of phase solution}.
\end{theorem}
Moreover, one can state a convergence result with linear speed, similar to Proposition~\ref{prop:linearcv}, for the Damped Newton algorithm.


\subsection{Point source far field}
Suppose rays emanate from the origin $O$, $\Pi$ is the plane $x_3=a$ and ${\bf m}$ is a unit direction. 
Then the metasurface $(\Pi,\phi)$ refracting rays from $O$ into the direction ${\bf m}$ (with $\nabla \phi(x)\cdot e=0$, $e=(0,0,1)$, i.e., 
$\phi$ tangential to $\Pi$\footnote{Each $\phi(x)=|x|-\m\cdot x+h(x_3)$ satisfies \eqref{eq:generalized snell law} with $n_1=n_2=1$ but is not tangential to $\Pi$ unless $h$ is constant.}) is given by
\begin{equation}\label{eq:far field phase}
\phi(x)=|x|-{\bf m}\cdot x+C
\end{equation}
where $x=(x_1,x_2,x_3)$ and $C$ a constant, see \cite[Section 4.A]{gps}.
Let $\Omega_0,\Omega'\subset S^2_+$ and let $\Omega=\{\lambda x:x\in \Omega_0, \lambda x\in \Pi\}$.

\begin{definition}[Admissible phase for the far field]\label{def:admissible phase far field}
The function $\phi:\Omega\to \R$ is a far field admissible phase refracting $\Omega$ into $\Omega'$ if for each $X_0\in \Omega$ there exists $b\in \R$ and $y\in \Omega'$ such that 
\[
\phi(X)\geq |X|-y\cdot X+b \quad \forall X\in \Omega, \quad \phi(X_0)= |X_0|-y\cdot X_0+b.
\]
When this happens, we say that $|X|-y\cdot X+b$ supports $\phi$ at $X_0$.
\end{definition}

In this case, the analysis about existence and uniqueness of solutions follows the lines already described and therefore we omit more details. It yields a theorem similar to Theorem \ref{thm:existencesumofdiracmassessecondbdryMA} now in terms of the far field supporting phases $|X|-y\cdot X+b$. 


We complete the paper mentioning that the phases for near field problem given by Definition \ref{def:admissible phase near field}
converge to the phases for far field problem in Definition \ref{def:admissible phase far field} when the target $T$ goes to infinity along fixed directions as indicated in the following lemma. 

\begin{lemma}\label{lm:convergence of costs}
We have the following convergence
\[
|X|+|X-P|-\(|P|+b\)\to |X|-\m\cdot X-b
\]
with $\m=P/|P|$ as $|P|\to \infty$ uniformly for $X\in K\subset \R^3$ compact and $b\in I$ a bounded interval.
\end{lemma}

\begin{proof}
Set $\Delta=|X-P|+|P|$ and write
\begin{align*}
|X-P|-\(|P|+b\)
&=
\dfrac{|X-P|^2-\(|P|+b\)^2}{\Delta+b}
=
\dfrac{|X|^2-2\(X\cdot P+b\,|P|\)-b^2}{\Delta+b}\\
&=
\dfrac{|X|^2}{\Delta+b}-\dfrac{2\,|P|\(X\cdot \m+b\)}{\Delta+b}-\dfrac{b^2}{\Delta+b}.
\end{align*}
Since $\Delta=\left|X-|P|\,\m\right|+|P|=|P|\(1+\left|\dfrac{X}{|P|}-\m \right| \)$, 
we have 
\[
\dfrac{\Delta +b}{2|P|}=
\dfrac{1+\left|\dfrac{X}{|P|}-\m \right|}{2}+\dfrac{b}{2|P|}\to 1
\]
as $|P|\to \infty$, uniformly for $X$ in a compact set and $b$ in a bounded interval, the lemma follows.
\end{proof}